\documentclass[journal,twocolumn,twoside]{ieeetran}
\ifCLASSINFOpdf
\else
\fi
	\input epsf  
	\usepackage[comma,numbers,square,sort&compress]{natbib}
	\usepackage{algorithm} 
	\usepackage{algorithmic} 
	\usepackage{setspace}
	\usepackage{arydshln}
	\usepackage{multirow} 
	\usepackage{xcolor}
	\usepackage{setspace}
	\usepackage{amsmath}
	
	\newtheorem{theorem}{\textbf{Theorem}}
	\newtheorem{lemma}{\textbf{Lemma}}
	\newtheorem{example}{\textbf{Example}}

	\newtheorem{remark}{\textbf{Remark}}
	\newtheorem{definition}{\textbf{Definition}}
	
	\newtheorem{proposition}{\textbf{Proposition}}
	
	\newcommand\myeq{\mathrel{\stackrel{\makebox[0pt]{\mbox{\normalfont\tiny (b)}}}{=}}}
	\newcommand\mge{\mathrel{\stackrel{\makebox[0pt]{\mbox{\normalfont\tiny (a)}}}{\ge}}}
	\newcommand\myeqa{\mathrel{\stackrel{\makebox[0pt]{\mbox{\normalfont\tiny (a)}}}{=}}}
	
	\usepackage{pmat}
		\newenvironment{proof}{{\noindent{\bf \emph{Proof:}}}\quad}{\hfill $\square$\par}
		
		\usepackage{multirow}
		\usepackage{url}
		\usepackage{subfigure}
		\usepackage{fancyhdr}
		\usepackage{amsmath}
		\usepackage{multirow}
		\usepackage{amssymb }
		\usepackage{color}
		\usepackage{graphics} 
		\usepackage{graphicx}
		\usepackage{geometry}
		\usepackage{setspace}

		\renewcommand{\baselinestretch}{1}

\begin{document}
			%
				%
				
				%
				\title{\LARGE {Perturbation-Tolerant Structural Controllability for Linear Systems}} 
				\author{ Yuan Zhang, Yuanqing Xia, Gang Wang, and Jinhui Zhang
					\thanks{ This work was supported in part by the  China Postdoctoral Innovative Talent Support Program under Grant BX20200055,  the China Postdoctoral
						Science Foundation under Grant 2020M680016, the National Natural Science Foundation of China under Grant 62003042, and the State Key Program of National Natural Science Foundation of China under Grant 61836001.  Yuan Zhang, Yuanqing Xia (corresponding author), Gang Wang and Jinhui Zhang are with the School of Automation, Beijing Institute of Technology, Beijing, China {(email: {\tt\small zhangyuan14@bit.edu.cn, xia\_yuanqing@bit.edu.cn, gangwang@bit.edu.cn, and	zhangjinh@bit.edu.cn}).} 
						
				} }
			\pagestyle{empty} 
			\maketitle
			\thispagestyle{empty} 
			{\small
				\begin{abstract}

					This paper proposes a novel notion termed perturbation-tolerant structural controllability (PTSC) to {study the generic property of controllability
					preservation for a structured linear system under structured perturbations}. More precisely, a structured system is said to be PTSC with respect to a given perturbation structure, if for almost all of its controllable realizations, there are no addable complex-valued perturbations with their zero/nonzero patterns prescribed by the perturbation structure that can make the resulting system uncontrollable.  We prove that genericity exists in this notion in the sense that, for almost all controllable realizations of a structured system, either there exist such addable structured perturbations rendering the resulting systems uncontrollable, or there is not such a perturbation. We give a decomposition-based necessary and sufficient condition for a single-input linear system to ensure PTSC, whose verification has polynomial time complexity. We then present some intuitive graph-theoretic conditions for PTSC. As an application, our results can serve as some feasibility conditions for the conventional {\emph{structured controllability radius}} problems from a generic view.
					
					
					
				\end{abstract}
				\begin{IEEEkeywords}
					Structural controllability, structured perturbations, controllability preservation, generic property
				\end{IEEEkeywords}
				\section{Introduction}
				In recent years, security has been becoming an attractive issue in the control and estimation of cyber-physical systems, such
				as chemical processes, power grids and transportation networks \cite{Buldyrev2009Catastrophic,fawzi2014secure,mitra2019byzantine}. The robustness of various system properties has been investigated under internal faults (like disconnections of links/nodes \cite{Buldyrev2009Catastrophic,zhang2020generic}) or external attacks (like adversarial sensor/actuator attacks \cite{fawzi2014secure}), including stability \cite{FabioFragility2018}, stabilization \cite{De2015Input}, controllability and observability \cite{commault2008observability,Rahimian2013Structural,zhang2019minimal}. Particularly, as a fundamental system property, controllability/observability under {\emph{structural perturbations}} (i.e., perturbations that can change the zero/nonzero patterns of the original system matrices, such as link/node/actuator/sensor removals or additions) has been extensively explored on its robust/resilence performance. To name a few, \cite{commault2008observability} considered observability preservation under sensor removals, \cite{Rahimian2013Structural} investigated controllability preservation under simultaneous link and node failures, while \cite{zhang2019minimal} systematically studied the involved optimization problems with respect to link/node/actuator/sensor removals from a computational perspective. {Since controllability/observability is a generic property in the sense that either the system realization is controllable for almost all values of the free parameters, or is not controllable for any values of the free parameters \cite{generic}, its robustness is overwhelmingly dominated by the robustness property of the graph representations associated with the system structure (see the next sections \ref{sec-II-B} and \ref{sec-III-C} for details). } 
				

				In the context of network systems, structural perturbation is a class of perturbations that make the corresponding links have zero weights. In the more general case where the perturbed links do not necessarily result in zero weights, controllability robustness has also been studied by computing the distance (in terms of the $2$-norm or the Frobenius norm) from a controllable system to the set of uncontrollable systems \cite{paige1981properties,WickDistance1991,Hu2004Real}. Such a notion, also named {\emph{controllability radius}}, was seemingly first proposed in \cite{paige1981properties}, and then developed by several other researchers on its efficient computations \cite{WickDistance1991,Hu2004Real}. Recently, by restricting the perturbation matrices to a prescribed structure, the so-called {\emph{structured controllability radius
						problem}} (SCRP), i.e., determining the smallest (Frobenius or $2$-) norm additive perturbation with a prescribed structure for which controllability fails to hold, has also attracted researchers' interest \cite{khare2012computing,johnson2018structured,bianchin2016observability,zhang2022real}. Towards this problem, various numerical algorithms have been proposed \cite{khare2012computing,johnson2018structured,bianchin2016observability}. However, due to the nonconvexity, most algorithms are suboptimal and iterative \cite{johnson2018structured}. Moreover, since most of them adopted some relaxation techniques in the iterating, and owing to the involved rounding errors and case-sensitive termination thresholds, there may even not be guarantees in practice that the returned perturbations can make the original system uncontrollable \cite{johnson2018structured}.
				
				On the other hand, to avoid the potential numerical issues, strong structural controllability (SSC), a notion proposed by Mayeda and Yamada \cite{mayeda1979strong}, could also be used to measure the controllability robustness of a system against numerical perturbations. In the SSC theory, the system parameters are divided into two categories, namely, indeterminate parameters and fixed zeros. A system is SSC, if whatever values (other than zero) the indeterminate parameters may take, the system is controllable. Criteria for SSC in the single-input case was given in \cite{mayeda1979strong}, and then extended to the multi-input cases in \cite{bowden2012strong,monshizadeh2014zero} and undirected networks in \cite{mousavi2017structural}. \cite{popli2019selective,jia2020unifying} reinvestigated SSC allowing the existence of parameters \! that \! can \! take \! arbitrary \! values including zero and nonzero. Note for the SSC theory to measure controllability robustness, the perturbations should have the same zero/nonzero structure as the original systems. While in practice, perturbations could happen only on partial system components (such as a subset of links of a network) and do not necessarily have the same structure as the original systems.
				
			 Motivated by the above observations, in this paper, we propose a novel notion to study controllability robustness under structured perturbations, namely, the perturbation-tolerant structural controllability (PTSC). A structured system is PTSC with respect to (w.r.t) a predefined perturbation structure, if almost all its controllable realizations can preserve controllability under arbitrary complex-valued perturbations with their zero/nonzero structure prescribed by the perturbation structure. The main contributions are as follows:

				1) We propose a novel notion of PTSC to study controllability preservation for a structured system under structured perturbations. For the first time, the perturbation structure in our notion can be {\emph{arbitrary}} w.r.t. the structure of the original system. This provides a new view in studying the robustness of structural controllability other than {\emph{structural perturbations}}, and has direct application in the feasibility of SCRPs (Sections \ref{ProblemFormulation} and \ref{apply-scrp}).
				
				
				2) We prove that genericity exists in this notion, in the sense that, for almost all controllable realizations of a structured system, either there exists an addable perturbation with the prescribed zero/nonzero structure rendering the resulting systems uncontrollable, or there is not such a perturbation (Section \ref{sec_gener}).

				
				3) We give a decomposition-based necessary and sufficient condition for a single-input system to be PTSC, whose verification has polynomial time complexity. The derivation is based on the   one-edge preservation principle and a series of nontrivial results on the roots of determinants of generic matrix pencils (Section \ref{sec-single}).
				
				4) We provide some intuitive graph-theoretic necessary and/or sufficient conditions for PTSC in the single-input case (Section \ref{graph-theoretic}).
				
				
				We believe the PTSC theory developed herein could be complementary to the existing SSC theory and provide a practical tool in estimating controllability robustness under structured perturbations. It is also beneficial in identifying critical/vulnerable edges in networks, i.e., edges whose change in weight could destroy system controllability. A conference version of this paper has appeared in \cite{zhang2021ptsc}, which is specific to the notion of PTSC for single-input systems only, and no graph-theoretic conditions are given therein. This paper provides a complete definition of PTSC both for single and multiple input systems, proving the genericity and providing proofs of results in \cite{zhang2021ptsc}.
				
				{\bf Outline}: Section II introduces the PTSC notion. Section III presents some preliminaries required for the further derivations. Section IV establishes the genericity within PTSC. Section V gives a necessary and sufficient condition for a single-input system to be PTSC, while Section VI presents some graph-theoretic conditions. Section VII discusses the application of PTSC to the SCRP. The last section concludes this paper.

				{\bf Notations}: Given a $p\in {\mathbb N}$, define ${\cal J}_p\doteq \{1,...,p\}$. For $q,p\in {\mathbb N}$ and $p\ge q$, define the set ${\cal J}_p^q\doteq \{{\cal I}\subseteq {\cal J}_p: |{\cal I}|=q\}$. For a $p\times q$ matrix $M$, $M[{\cal I}_1,{\cal I}_2]$ denotes the submatrix of $M$ whose rows are indexed by ${\cal I}_1$ and columns by ${\cal I}_2$, ${\cal I}_1\subseteq {\cal J}_p$, ${\cal I}_2\subseteq {\cal J}_q$. $||M||_0$ represents the total number of nonzero entries in $M$. $\det (\cdot)$ stands for the determinant.   For a vector $v$, $v_i$ denotes the $i$th entry.  

				
				\section{Problem Formulation} \label{ProblemFormulation}
				\subsection{Structured Matrix} 
				{A structured matrix is a matrix whose entries are either fixed zero or not fixed to zero, and is often represented by a matrix with entries chosen from $\{0,*\}$. Here, $0$ stands for a fixed zero entry, $*$ for an entry that is not fixed to zero. The latter is called a nonzero entry in the sequel.} 
				 Let $\{0,*\}^{p\times q}$ be the set of all $p\times q$ dimensional structured matrices. For $\bar M \in \{0,*\}^{p\times q}$, define the set of matrices 
				\begin{equation}\label{structure-matrix1}{\bf S}_{\bar M}=\left\{ M\in {\mathbb C}^{p\times q}: M_{ij}=0\ {\text{if}}\ {\bar M}_{ij}=0 \right\},\end{equation}
				Any $M\in {\bf S}_{\bar M}$ is called a {\emph{realization}} of $\bar M$.  	A generic realization of $\bar M$ is a realization whose $*$ entries are assigned with free parameters (i.e., parameters that can take values independently). We also call such a matrix a {\emph{generic matrix}} \cite{Murota_Book} without specifying the corresponding structured matrix. For a generic matrix $M$ and a constant matrix $N$ with the same dimension, $M+\lambda N$ defines a generic matrix pencil, which can be seen as a {\emph{matrix-valued polynomial of free parameters in $M$ and the variable $\lambda$}}. 
				For two structured matrices $\bar M, \bar N \in \{0,*\}^{p\times q}$, $\vee$ is the entry-wise OR operation, i.e., $(\bar M\vee \bar N)_{ij}= *$ if $\bar M_{ij}=*$ or $\bar N_{ij}=*$, and $(\bar M\vee \bar N)_{ij}=0$ otherwise.  With some abuse of notation, we say $\bar M\subseteq \bar N$, if $\bar M_{ij}=*$ implies that $\bar N_{ij}= *$.

				\subsection{Notion of PTSC} \label{sec-II-B}
				Consider the linear time invariant  system
				\begin{equation}\label{plant}\dot x(t)=Ax(t)+Bu(t),\end{equation}
				where $x\in {\mathbb C}^n$ is the state vector, $u\in {\mathbb C}^m$ is the input vector, $A\in {\mathbb C}^{n\times n}$, and $B\in {\mathbb C}^{n\times m}$. For notation simplicity, we may use $(A,B)$, or alternatively $[A,B]$ (if not causing confusion), to represent system (\ref{plant}). It is known that $(A,B)$ is controllable (i.e., system (\ref{plant}) is controllable), if and only if the controllability matrix ${\cal C}(A, B)$ defined as follows has full row rank
				$${\cal C}(A,B)=[B,AB,\cdots, A^{n-1}B].$$
				
				Let ${\bar F}\in \{0,*\}^{n\times {(n+m)}}$ be a structured matrix specifying the structure (i.e., zero/nonzero pattern)
				of the perturbation (matrix) $[\Delta A, \Delta B]$, that is, ${\bar F}_{ij}=0$ implies $[\Delta A, \Delta B]_{ij}=0$. In other words, $[\Delta A, \Delta B]\in {\bf S}_{\bar F}$, with ${\bf S}_{\bar F}$ defined in (\ref{structure-matrix1}). We may occasionally use $[\Delta \bar A, \Delta \bar B]$ to denote $\bar F$ and also write $(\Delta A, \Delta B)\in {\bf S}_{\bar F}$. It is emphasized that in this paper, the perturbations are allowed to be {\emph {complex-valued}. We will briefly discuss the case where the perturbations are restricted to be in the real field in Section \ref{sec_gener}. {\emph{ We call $(A,B)$ the original system, and $(A+\Delta A, B+\Delta B)$ the perturbed system.  }}					Let $\bar A \in \{0,*\}^{n\times n}, \bar B\in \{0,*\}^{n\times m}$ be the structured matrices specifying the zero/nonzero pattern of $A,B$. Put differently, $A\in {\bf S}_{
						\bar A}$ and $B\in  {\bf S}_{\bar B}$.

					\begin{definition}[PTC]
						System $(A,B)$ in (\ref{plant}) is said to be {\emph{perturbation-tolerantly controllable}} (PTC) w.r.t. $\bar F$, if for all $[\Delta A,\Delta B]\in {\bf S}_{\bar F}$, $(A+\Delta A, B+\Delta B)$ is controllable. If $(A,B)$ is controllable but not PTC w.r.t. $\bar F$ (i.e., there exists a $[\Delta A, \Delta B]\in {\bf S}_{\bar F}$ making $(A+\Delta A, B+\Delta B)$ uncontrollable), $(A,B)$ is said to be {\emph{perturbation-sensitively controllable}} (PSC) w.r.t. $\bar F$.
					\end{definition}

					\begin{definition}[Structural controllability] $(\bar A, \bar B)$ is said to be structurally controllable, if there exists a realization $[A,B]\in {\bf S}_{[\bar A, \bar B]}$ so that $(A,B)$ is controllable.
					\end{definition}
					
					A property is called generic, if either for almost all\footnote{Hereafter, by `almost all' we mean `all except a set of Lebesgue measure zero in the corresponding parameter space'. In other words, let $\cal S$ contain the parameter values that make the considered property hold. Then, ${\cal S}$ can be written as ${\mathbb F}^d\backslash {\cal P}$ for some proper algebraic variety of the whole parameter space ${\mathbb F}^d$, where $d$ is the number of parameters in the system parametrization (a proper algebraic variety is the zero set of some nontrivial polynomial, which has Lebesgue measure zero \cite{dummit2004abstract}).} parameter values, this property holds, or for almost all parameter values, this property does not hold \cite{generic}. It is well-known that controllability is a generic property in the sense that, if $(\bar A, \bar B)$ is structurally controllable, then almost all realizations of $(\bar A, \bar B)$ are controllable; otherwise, none realization of $(\bar A, \bar B)$ is controllable (see \citep[Prop 1]{C.T.1974Structural} or \citep[Prop 1]{zhang2019structural} for argument of this assertion). For a structurally controllable pair $(\bar A, \bar B)$, let ${\bf CS}(\bar A,\bar B)$ denote the set of all {\emph{controllable}} complex-valued realizations of $(\bar A, \bar B)$. With the notions above, we give the definition of PTSC as follows.
					

					\begin{definition}[PTSC]\label{DefinitionForPTSC} Given $(\bar A, \bar B)$ and a perturbation structure ${\bar F}$, $(\bar A, \bar B)$ is said to be PTSC w.r.t. $\bar F$, if for almost all $(A,B)\in {\bf CS}(\bar A,\bar B)$, $(A,B)$ is PTC w.r.t. $\bar F$. If $(\bar A,\bar B)$ is structurally controllable but not PTSC w.r.t. $\bar F$, $(\bar A,\bar B)$ is alternatively called perturbation-sensitively structurally controllable (PSSC) w.r.t. $\bar F$.
					\end{definition}

				{\emph{In the rest of this paper, for ease of description, when referring to PTSC or PSSC, we implicitly assume that the original (unperturbed) system is structurally controllable.}}	In Section \ref{sec_gener}, we will prove that PTC (as well as PSC) is a generic property, that is, depending on the structure of $(\bar A, \bar B)$ and $\bar F$, either for almost all $(A,B)\in {\bf CS}(\bar A, \bar B)$, $(A,B)$ is PTC w.r.t. $\bar F$, or for almost all $(A,B)\in {\bf CS}(\bar A, \bar B)$, $(A,B)$ is PSC w.r.t. $\bar F$. Hence, the definition of PTSC makes sense. Before doing so, we  give an illustrative example of PTSC as well as PSSC.

					\begin{example}\label{exp1}Consider a single-input system $(A,b)$ parameterized by $(a_{21},a_{32},a_{42},a_{44},a_{41},b_1)\in {\mathbb C}^6$ as
						{\small$$
							A=\left[
							\begin{array}{cccc}
								0 & 0 & 0 & 0 \\
								a_{21} & 0 & 0 & 0 \\
								0 & a_{32} & 0 & 0 \\
								a_{41} & a_{42} & 0 & a_{44} \\
							\end{array}
							\right],b=\left[
							\begin{array}{c}
								b_1 \\
								0 \\
								0 \\
								0 \\
							\end{array}
							\right].
							$$}Two perturbations $[\Delta A_i, \Delta b_i]$ ($i=1,2$) are given as {\small
							$$\left[
							\begin{array}{ccccc}
								0 & 0 & {\Delta a_{13}} & {\Delta a_{14}} & 0 \\
								0 & 0 & 0 & 0 & 0 \\
								0 & 0 & 0 & 0 & 0 \\
								0 & 0 & 0 & 0 & 0 \\
							\end{array}
							\right],\left[
							\begin{array}{ccccc}
								0 & 0 & 0 & 0 & 0 \\
								0 & 0 & 0 & 0 & 0 \\
								0 & 0 & {\Delta a_{33}} & 0 & 0 \\
								0 & 0 & 0 & 0 & {\Delta b_{4}} \\
							\end{array}
							\right],$$}where ${\Delta a_{13}},{\Delta a_{14}},{\Delta a_{33}},{\Delta b_{4}}$ can be arbitrarily selected. It can be obtained $\det {\cal C}(A+\Delta A_1, b+\Delta b_1)=a_{21}^2a_{32}b_1^4(a_{41}a_{44}^2 + a_{21}a_{42}a_{44})$, which is independent of $\Delta a_{13}$ and $\Delta a_{14}$. Hence, whatever values $[\Delta A_1, \Delta b_1]$ may take, the perturbed system is always controllable (provided $(A,b)$ is controllable). By definition, the structured system  $(\bar A, \bar b)$ is PTSC w.r.t. the structured perturbation $[ \Delta \bar A_1, \Delta \bar  b_1]$. On the other hand, $\det {\cal C}(A+\Delta A_2, b+\Delta b_2)=a_{21}^2a_{32}b_{1}^3(a_{44}^3\Delta b_{4} - a_{44}^2\Delta b_{4}\Delta a_{33} + a_{44}^2a_{41}b_{1} - a_{44}a_{41}b_{1}\Delta a_{33} + a_{21}a_{42}a_{44}b_{1} - a_{21}a_{42}b_{1}\Delta a_{33})$. It can be seen, in case $a_{44}^3 - \Delta a_{33}a_{44}^2\ne 0$ (or $\Delta b_{4}a_{44}^2+a_{41}b_{1}a_{44}+a_{21}a_{42}b_{1}\ne 0$) and $(A,b)$ is controllable, there is a value for ${\Delta b_{4}}$ (resp. ${\Delta a_{33}}$) making $\det {\cal C}(A+\Delta A_2, b+\Delta b_2)=0$, rendering the perturbed system uncontrollable. This indicates, the structured system $(\bar A, \bar b)$ is PSSC w.r.t. the perturbation $[ \Delta \bar A_2,  \Delta \bar b_2]$. 
					\end{example}
				
				    \begin{remark}If $\bar F$ has the same zero/nonzero pattern as $[\bar A, \bar B]$, then $(\bar A, \bar B)$ is automatically PSSC w.r.t. $\bar F$, as in this case all $*$ entries in $[\bar A, \bar B]$ can be perturbed to zero. In the context of network systems, it often holds $\bar F\subseteq [\bar A, \bar B]$, corresponding to that a partial set of edges have their weights perturbed while the rest remain unchanged. 
				   	\end{remark}
					
					\begin{remark}[Alternative definition of PTSC]\label{second-viewpoint}
						We may also revisit PTSC from the standpoint of the {\emph{perturbed structured system}}. In this system, entries of system matrices can be divided into three categories, namely, the {\emph{fixed zero entries}}, the {\emph{unknown generic entries}} ($*$ entries in $[\bar A, \bar B]$) that take fixed but unknown values (they can be modeled as algebraically independent parameters), and the {\emph{perturbed entries}} ($*$ entries in $\bar F$) which can take arbitrarily complex values.\footnote{An unknown generic entry and a perturbed one are allowed to be the same entry. In this case, it suffices to regard this entry as a perturbed one.}  PTSC of this system requires that for almost all values of the unknown generic entries making the original system controllable, the corresponding perturbed systems are controllable for arbitrary values of the perturbed entries. 
					\end{remark}
				
				\begin{remark} \label{justification}
					It is worth mentioning that assuming algebraic independence among the considered parameters to study certain generic properties related to the problem structure is common in the 
					mathematics community \cite{dummit2004abstract}. A similar notion to the PTSC in spirit is the widely-accepted notion generic low-rank matrix completion (GLRMC) \cite{kiraly2015algebraic,zhang2021generic}. GLRMC is typically formulated as, given a matrix with a known pattern of observed and missing entries, verifying whether for almost all values of the observed entries, there is a matrix completion (obtained by assigning values to the missing entries) with a predefined rank. Note the problem of PTSC differs from the GLRMC, as the former essentially reduces to the problem of completing a {\emph{matrix pencil}} $M+\lambda N$ (rather than a matrix) such that its minimum rank over $\lambda \in {\mathbb C}$ is generically less than $n$ ($n$ is the number of rows of $M$). We shall discuss one application of PTSC in Section \ref{apply-scrp}. 
				\end{remark}
				
				In the sequel, we may directly say a system is PTSC, when either the corresponding perturbation structure is clear from the context, or this system refers to a perturbed system with three categories of entries as in Remark \ref{second-viewpoint}. The main purposes of this paper are to i) justify the genericity of PTC; ii) provide necessary or sufficient conditions for PTSC; and iii) show the applications of PTSC in SCRPs.
					

					\subsection{Relations with SSC} \label{sec-scc}
					We shall briefly discuss the relation between SSC and PTSC here.
					\begin{definition}[SSC,\cite{mayeda1979strong}] $(\bar A, \bar B)$ is said to be SSC, if every realization of $(\bar A, \bar B)$ is controllable subject to that each $*$ entry takes strictly nonzero values. 
					\end{definition}
					
					As mentioned earlier, SSC could be seen as the ability of a system to preserve controllability under perturbations that have the same structure as the system itself, with the constraint that the perturbed entries of the resulting system cannot be zero. Combined with Remark \ref{second-viewpoint},  the essential difference between PTSC and SSC lies in two aspects:
					\begin{itemize}
						\item[i)] The perturbed entries can take arbitrary values including zero in PTSC, while they must take nonzero values in SSC;
						\item[ii)] In SSC, all nonzero entries should be perturbed, while in PTSC, an arbitrary subset of entries (prescribed by the perturbation structure) can be perturbed and the remaining entries remain unchanged.
					\end{itemize}Because of them, neither criteria for SSC can be converted to those for PTSC, nor the converse. Recently, some researchers have extended the SSC to allowing entries that can take arbitrary values \cite{jia2020unifying, popli2019selective}. Built on such an extension, item ii) identified above remains an essential difference between SSC and PTSC. 
				
				\begin{example}
					To show the SSC in \cite{jia2020unifying, popli2019selective} does not subsume PTSC, we still consider $(A,b)$ and the perturbation $(\Delta A_1, \Delta b_1)$ in Example \ref{exp1}, with their corresponding structured matrices $(\bar A, \bar b)$ and $(\Delta \bar A_1, \Delta \bar b_1)$. Let the $*$ entries in $[\Delta \bar A_1, \Delta \bar b_1]$ take
					arbitrary values, while in $[\bar A, \bar b]$ take nonzero values. Because there are nonzero values for $a_{41}, a_{44}, a_{21}, a_{42}$ making $\det {\cal C}(A+\Delta A_1, b+\Delta b_1)=0$, we know $[\bar A, \bar b]\vee [\Delta \bar A_1, \Delta \bar b_1]$ is not SSC in the sense of \cite{jia2020unifying, popli2019selective}. By contrast, Example \ref{exp1} shows 
					$[\bar A, \bar b]\vee [\Delta \bar A_1, \Delta \bar b_1]$ is indeed PTSC. 
				\end{example}
					


%

					\section{Preliminaries and Terminologies}
					In this section, we introduce some preliminaries as well as terminologies in graph theory and structural controllability.
					\vspace{-0.4cm}
					\subsection{Graph Theory}
					If not specified, all graphs in this paper refer to directed graphs. A graph is denoted by ${\cal G}=({\cal V},{\cal E})$, where ${\cal V}$ is the vertex set, and ${\cal E}\subseteq {\cal V}\times {\cal V}$ is the edge set. Let ${\cal E}_{x}^{\rm in}$ and ${\cal E}_{x}^{\rm out}$ be the set of ingoing edges into $x\in {\cal V}$ and outgoing edges from $x$, respectively, i.e., ${\cal E}^{\rm in}_{x}=\{z\in {\cal V}: (z,x)\in {\cal E}\}$, ${\cal E}^{\rm out}_{x}=\{z\in {\cal V}: (x,z)\in {\cal E}\}$.  A {\emph{path}} from vertex $v_i$ to vertex $v_j$ is a sequence of distinct edges $(v_i,v_{i+1})$, $(v_{i+1},v_{i+2})$, $\cdots$, $(v_{j-1},v_{j})$ where each edge belongs to $\cal E$. A path from a vertex to itself is called a {\emph{cycle}}. A path not containing any cycle is an {\emph{elementary path}}. A strongly connected component (SCC) of $\cal G$ is a subgraph so that there is a path from every of its vertex to the other, and no more vertex can be included in this subgraph without breaking such property. We may also call the vertex set of this subgraph an SCC. An {\emph{induced subgraph}} of $\cal G$ by ${\cal V}_s\subseteq {\cal V}$ is a graph formed by ${\cal V}_s$ and all of the edges in $({\cal V}_s\times {\cal V}_s) \cap {\cal E}$. For a set ${\cal E}_s\subseteq {\cal E}$, ${\cal G}-{{\cal E}_s}$ denotes the graph obtained from $\cal G$ after deleting the edges in ${\cal E}_s$; similarly, for ${\cal V}_s\subseteq {\cal V}$, ${\cal G}-{\cal V}_s$ denotes the graph after deleting vertices in ${\cal V}_s$ and all edges incident thereto. For two graphs ${\cal G}_i=({\cal V}_i, {\cal E}_i)$, $i=1,2$, ${\cal G}_1\cup {\cal G}_2$ denotes the graph $({\cal V}_1\cup {\cal V}_2, {\cal E}_1\cup {\cal E}_2)$. 
					

						A bipartite graph is denoted by $({\cal V}_1, {\cal V}_2, {\cal E})$, with vertex bipartitions ${\cal V}_1, {\cal V}_2$ and edge set ${\cal E}$.	
				 A matching of a bipartite graph is a subset of its edges among which any two do not share a common vertex. The{\emph{ maximum matching}} is the matching containing as many edges as possible. The number of edges contained in a maximum matching of $\cal G$ is denoted by ${\rm mt}({\cal G})$. For a weighted bipartite graph $\cal G$, where each edge is assigned a non-negative weight, the weight of a matching is the sum of all its edge weights. Denote by ${\rm minWMM({\cal G})}$ (resp. ${\rm maxWMM({\cal G})}$) the minimum (resp. maximum) weight over all maximum matchings of ${\cal G}$. The bipartite graph associated with a matrix $M$ is given by ${\mathcal B}(M)=({\cal V}^+ ,{\cal V}^- ,{\cal E})$, where the left (resp. right) vertex set ${\cal V}^+$ (resp. ${\cal V}^-$) corresponds to row (resp. column) set of $M$, and the edge set corresponds to the set of nonzero entries of $M$, i.e., ${\mathcal E}=\{(i,j): i\in {\cal V}^+, j\in {\cal V}^-, M_{ij}\ne 0\}$.

					The {\emph{generic rank}} of a structured matrix $\bar M$, given by ${\rm grank}(\bar M)$, is the maximum rank $\bar M$ can achieve as the function of its free parameters, which also equals the rank it achieves for almost all choices of its parameter values \citep[Page 38]{Murota_Book}. 
					It is known that ${\rm grank}(\bar M)$ equals the cardinality of the maximum matching of ${\mathcal B}(\bar M)$ \citep[Props 2.1.12 and 2.2.25]{Murota_Book}. For a polynomial-valued matrix $M$, its generic rank ${\rm grank}(M)$ is defined as the maximum rank it can achieve as the function of its free parameters.

					{\emph{Dulmage-Mendelsohn decomposition}} (DM-decomposition) is a decomposition of a bipartite graph w.r.t. maximum matchings. Let
					${\cal G}=({\cal V}^{+}, {\cal V}^{-}, {\cal E})$ be a bipartite graph.  An edge is said to be admissible,
					if it is contained in some maximum matching of~$\cal G$. 
					
					\begin{definition}[DM-decomposition, \cite{Murota_Book}] \label{DM-def} The DM-decomposition of a bipartite graph
						${\cal G}=({\cal V}^{+}, {\cal V}^{-}, {\cal E})$ is to decompose $\cal G$ into subgraphs ${\cal G}_k=({\cal V}^+_k,{\cal V}^-_k,{\cal E}_k)$ ($k=0,1,...,d,\infty$)
						(called DM-components of $\cal G$) satisfying:	1) ${\cal V}^{\star}=\bigcup \nolimits_{k=0}^{\infty} {\cal V}^{\star}_k$, ${\cal V}^{\star}_i\bigcap {\cal V}^{\star}_j=\emptyset$ for $i\ne j$,  with $\star= +$, $-$; ${\cal E}_k=\{(v^+,v^-)\in {\cal E}: v^+\in {\cal V}^+_k, v^-\in {\cal V}^-_k\}$;
		        	2) For $1\le k \le d$ (consistent components): ${\rm mt}({\cal G}_k)=|{\cal V}_k^+|=|{\cal V}_k^-|$, and each $e\in {\cal E}_k$ is admissible in ${\cal G}_k$;
						for $k=0$ (horizontal tail): ${\rm mt}({\cal G}_0)=|{\cal V}_0^+|$, $|{\cal V}_0^+|<|{\cal V}_0^-|$ if ${\cal V}_0^+\ne \emptyset$, and each $e\in {\cal E}_0$ is
						admissible in ${\cal G}_0$; for $k=\infty$ (vertical tail): ${\rm mt}({\cal G}_{\infty})=|{\cal V}^-_\infty|$, $|{\cal V}_\infty^+|>|{\cal V}_\infty^-|$ if ${\cal V}_\infty^-\ne \emptyset$, and each $e\in {\cal E}_\infty$ is
						admissible in ${\cal G}_\infty$;
					   3) ${\cal E}_{kl}=\emptyset$ unless $1\le k\le l\le d$, and ${\cal E}_{kl}\ne \emptyset$ only if
						$1\le k\le l\le d$, where ${\cal E}_{kl}=\{(v^+,v^-)\in {\cal E}: v^+\in {\cal V}_k^+, v^-\in {\cal V}_l^-\}$;
						4) ${\cal M}$ ($\subseteq {\cal E}$) is a maximum matching of $\cal G$ if and only if ${\cal M}\subseteq \bigcup \nolimits_{k=0}^{\infty} {\cal E}_k$ and ${\cal M}\cap {\cal E}_k$ is a maximum matching of ${\cal G}_k$;
						5) ${\cal G}$ cannot be decomposed into more components satisfying 1)-4).
					\end{definition}
					
					For an $m\times l$ matrix $M$, the DM-decomposition of ${\cal B}(M)$ into graphs ${\cal G}_k=({\cal V}^+_k,{\cal V}^-_k,{\cal E}_k)$ ($k=0,1,...,d,\infty$) corresponds to that there exist two permutation matrices $P\in {\mathbb R}^{m\times m}$ and $Q\in {\mathbb R}^{l\times l}$ so that $PMQ$ is an upper block-triangular matrix, whose diagonal blocks $M_k=M[{\cal V}^+_k,{\cal V}^-_k]$ correspond to ${\cal G}_k$ ($k=0,1,...,d,\infty$).\footnote{For notation simplicity, we may use $M[{\cal V}_k^+,{\cal V}_k^-]$ to denote the submatrix of $M$ whose rows correspond to the vertex set ${\cal V}_k^+$ and columns to ${\cal V}_k^-$.} Such a process is also called the DM-decomposition of $M$. 
					
					A bipartite graph is said to be DM-irreducible if it cannot be decomposed into more than one nonempty component in the DM-decomposition.
					{\emph{DM-decomposition will play an important role in characterizing conditions for PTSC due to its close relation with the irreducibility of the determinant of a generic matrix}}. A multivariable polynomial $f$ is irreducible if it cannot be factored as $f=f_1f_2$ with $f_1,f_2$ being polynomials with smaller degrees than $f$.
					
					\begin{lemma} \citep[Theo 2.2.24, 2.2.28]{Murota_Book}\label{reduciblility} For a bipartite graph ${\cal B}(M)=({\cal V}^+ ,{\cal V}^- ,{\cal E})$ associated with a generic square matrix $M$, the following conditions are equivalent:	
						1) ${\cal B}(M)$ is DM-irreducible;
					    2) ${\rm mt}({\cal B}(M)-\{v_1,v_2\})= {\rm mt}({\cal B}(M))-1$ for any $v_1\in {\cal V}^+$ and $v_2\in {\cal V}^-$;
						3) $\det (M)$ is irreducible.
					\end{lemma}
					
					\subsection{Structural Controllability} \label{sec-III-C}
					It is convenient to map structured systems to graph elements. Given $(\bar A,\bar B)$, let $\mathcal{X}$, $\mathcal{U}$ denote the sets of state vertices and input vertices respectively, i.e., $\mathcal{X}=\{x_1,...,x_n\}$, $\mathcal{U}=\{{x_{n+1},...,x_{n+m}}\}$. Denote the edges by $\mathcal{E}_{\bar A}=\{(x_i,x_j): \bar A_{ji}\ne 0\}$, $\mathcal{E}_{\bar B}=\{({x_{n+j}},x_i): \bar B_{ij}\ne 0\}$.  Let $\mathcal{G}(\bar A,\bar B)= (\mathcal{X}\cup \mathcal{U}, \mathcal{E}_{\bar A}\cup \mathcal{E}_{\bar B})$ be the system graph of $(\bar A, \bar B)$. A state vertex $x\in\mathcal{X}$ is said to be input-reachable, if there exists a path from an input vertex to it in $\mathcal{G}(\bar A,\bar B)$. ${\cal G}({\bar F})=({\cal X}\cup {\cal U}, {\cal E}_{\bar F})$ is defined similarly, with ${\cal E}_{\bar F}=\{(x_i,x_j): {\bar F}_{ji}\ne 0\}$.
					
					In $\mathcal{G}(\bar A,\bar B)$, a {\emph{stem}} is an elementary path starting from an input vertex.  A bud is a cycle in ${\cal E}_{\bar A}$ with an {\emph{additional}} edge that ends (but not begins) in a vertex of this cycle. A {\emph{cactus}} is a subgraph ${\cal G}_s$ that is defined recursively as follows:  A cactus is a stem or is obtained by adding a bud to a smaller cactus ${\cal G}'_s$, in which the beginning vertex of this bud can be any vertex of ${\cal G}'_s$ except the last vertex of the stem contained in ${\cal G}'_s$.
					
					\begin{lemma}\citep[Theo 1]{generic} \label{theo-strucon} For system $(\bar A, \bar B)$ in (\ref{plant}), the following statements are equivalent:
						i) $(\bar A, \bar B)$ is structurally controllable;
							ii) ${\cal G}(\bar A, \bar B)$ is spanned by a collection of disjoint cacti;
						iii) Every state vertex is input-reachable, and ${\rm grank}([\bar A,\bar B])=n$.
					\end{lemma}
					
					\section{Genericity of PTC} \label{sec_gener}
					In this section, we prove the genericity of PTC. We also point out the subtle difference between the genericity of PTC in the single-input case and the multi-input one, and discuss PTC in the real field.  To simplify descriptions, hereafter, by $(A,B)$ we refer to a system that can be either single-input or multi-input, while by $(A,b)$ to a single-input system.

					
					\begin{theorem}[Genericity of PTC]\label{generic_PSC_multi} Given a structurally controllable pair $(\bar A, \bar B)$ and a perturbation structure $\bar F$, either for {\emph{almost all}} $(A,B)\in {\bf CS}(\bar A,\bar B)$, $(A,B)$ is PTC w.r.t. $\bar F$, or for {\emph{almost all}} $(A,B)\in {\bf CS}(\bar A,\bar B)$, $(A,B)$ is PSC w.r.t. $\bar F$.
					\end{theorem}
					
					The result above reveals that PTC (as well as PSC) is a generic property in ${\bf CS}(\bar A, \bar B)$.  The proof is postponed to Appendix A. Our proof adopts a similar technique to \cite{zhang2021generic}, which studies a different problem in low-rank matrix completions. In the following, we show for sinlge-input systems, PTSC has a stronger implication than Theorem \ref{generic_PSC_multi} in the sense that the first `almost all' therein could be changed to `all'.

					\begin{proposition}[Genericity in single-input case]\label{generic_PSC_single} With $(\bar A, \bar b)$ and $\bar F$ defined above, either for {\emph{all}} $(A,b)\in {\bf CS}(\bar A,\bar b)$, $(A, b)$ is PTC w.r.t. $\bar F$, or for {\emph{almost all}} $(A,b)\in {\bf CS}(\bar A,\bar b)$, $(A,b)$ is PSC w.r.t. $\bar F$.
					\end{proposition} 

					\begin{proof} Let $p_1,...,p_r$ be values that the $*$ entries of $[A,b]$ take, and $\bar p_1,...,\bar p_l$ be values that the perturbed entries of $[\Delta A, \Delta b]\in {\bf S}_{\bar F}$ take, with $r$, $l$ the numbers of the respective entries. Denote by $p\doteq (p_1,...,p_r)$ and $\bar p\doteq (\bar p_1,...,\bar p_l)$. It turns out $\det {\cal C}(A+\Delta A, b+\Delta b)$ can be factored as
						\begin{equation} \label{determinant_of_Ab}
							\det {\cal C}(A+\Delta A, b+\Delta b)=f(p)g(\bar p)h(p,\bar p),
						\end{equation}
						where $f(p)$ (resp. $g(\bar p)$) denotes the polynomial of $p$ (resp. $\bar p$) with real coefficients, and $h(p,\bar p)$ denotes the polynomial of $p$ and $\bar p$, in which at least one $p_i$ ($i\in \{1,...,r\}$) as well as one $\bar p_j$ ($j\in \{1,...,l\}$) has a degree no less than one.
						
						It can be seen that, if neither $g(\bar p)$ nor $h(p,\bar p)$ exists in the right-hand side of (\ref{determinant_of_Ab}), then for all $(A,b)\in {\bf CS}(\bar A,\bar b)$, $(A,b)$ is PTC w.r.t. $\bar F$, as in this case, (\ref{determinant_of_Ab}) is independent of $\bar p$. Otherwise, suppose there exists a $\bar p_j$, $j\in\{1,...,l\}$, so that there is a term $\bar p_j^k$ with degree $k\ge 1$ in $g(\bar p)h(p,\bar p)$. When $p$ and $\bar p \backslash \{\bar p_j\}$ take such values that the coefficient of $\bar p_j^k$, which is polynomial in $p$ and $\bar p \backslash \{\bar p_j\}$, is nonzero (referred to as (a1)), and meanwhile $(A,b)$ is controllable (referred to as (a2)), then according to the fundamental theorem of algebra (cf. \cite{dummit2004abstract}), there exists a complex value of $\bar p_j$ making $g(\bar p)h(p,\bar p)=0$, thus making $\det {\cal C}(A+\Delta A, b+\Delta b)=0$, leading to the uncontrollability of $(A+\Delta A, b+\Delta b)$. Note the set of values for $p$ and $\bar p\backslash \{\bar p_j\}$ satisfying (a1) has full dimension in ${\mathbb C}^{r+l-1}$. Thus, the set of values for $p$ validating (a1) forms a proper algebraic variety in ${\mathbb C}^{r}$, which thereby has zero Lebesgue measure in the set of values satisfying (a2). 
					\end{proof}

					
					Proposition \ref{generic_PSC_single} implies that, if there is a $(A,b)\in {\bf CS}(\bar A, \bar b)$ that is PSC w.r.t. $\bar F$, then almost all $(A,b)\in {\bf CS}(\bar A, \bar b)$ are PSC w.r.t. $\bar F$; otherwise, none $(A,b)\in {\bf CS}(\bar A, \bar b)$ is PSC w.r.t. $\bar F$. Compared with Theorem \ref{generic_PSC_multi}, this may not be surprising, since from the respective proofs, PTSC for single-input systems is directly related to the solvability of a single polynomial equation, while in the multi-input case, to the common zeros of a series of polynomial equations.
					
					\begin{example}[Single-input versus multi-input]
						Consider {\footnotesize$A=\left[
							\begin{array}{cc}
								a_{11} & a_{12} \\
								a_{21} & a_{22} \\
							\end{array}
							\right]
							$, $B=\left[
							\begin{array}{cc}
								b_{11} & 0 \\
								0 & b_{22} \\
							\end{array}
							\right]
							$}, and {\footnotesize$[\Delta A,\Delta B]=\left[
							\begin{array}{cccc}
								0 & \Delta a_{12} & 0 & 0 \\
								0 & 0 & 0 & 0 \\
							\end{array}
							\right]
							$}. If $b_{11}b_{22}\ne0$, it is clear that whatever $\Delta a_{12}$ takes, the perturbed system $(A+\Delta A,B+ \Delta B)$ is controllable. This means, $[\bar A, \bar B]$ is PTSC w.r.t. $[\Delta \bar A,\Delta \bar B]$.  On the other hand, when $b_{11}=0$ and $b_{22}\ne 0$, it can be validated that if $\Delta a_{12}=-a_{12}$, then $(A+\Delta A, B+\Delta B)$ is uncontrollable. Moreover, when $b_{11}\ne 0$ and $b_{22}=0$, $(A,B)$ collapses to a single-input system $(A,b)$. In this case,  $\det{\cal C}(A+\Delta A,b+\Delta b)=a_{21}b^2_{11}$, meaning that $(A,b)$ is PTC w.r.t. $[\Delta \bar A, \Delta \bar b]$, provided $(A,b)$ is controllable. This exemplifies, unlike the single-input case, even a multi-input structured system is PTSC, it cannot guarantee that all its controllable realizations are PTC w.r.t. the corresponding perturbations. 
					\end{example}
					
					Theorem \ref{generic_PSC_multi} indicates that it is the structure of the original system $(A,B)$ and the perturbation $\bar F$ that dominates the property of being PTC or PSC. This justifies the notion of PTSC in Definition \ref{DefinitionForPTSC}. Moreover, from the proof of Theorem \ref{generic_PSC_multi}, it is readily seen that even when the set ${\bf CS}(\bar A,\bar B)$ of the original system matrices is defined in the real field, provided that the perturbed entries can take complex values, genericity of PTC (PSC) still holds in the corresponding space. Further, the proof technique of Theorem \ref{generic_PSC_multi} can be trivially extended to show that, 
				    given that the unknown generic entries depend on some elementary parameters affinely, while there are some known algebraic dependence among the perturbed entries (such as symmetry), PTC remains a generic property in the space of the elementary parameters.
					
					
					However, if the original systems $(A,B)$ and the perturbed systems $(A+\Delta A, B+\Delta B)$ are both restricted to be in the real field, it is possible that for the same $(\bar A, \bar B)$ and $\bar F$, there is a subset of realizations of $(\bar A,\bar B)$ corresponding to a full-dimensional semi-algebraic subset of the real vector space\footnote{A semi-algebraic set is a subset of ${\mathbb R}^n$ for some real closed field defined by a finite sequence of polynomial equations and inequalities, or any finite union of such sets \cite{Sturmfels2002}.} that are PTC w.r.t. $\bar F$; meanwhile, the remaining realizations of $(\bar A,\bar B)$, also corresponding to a full-dimensional semi-algebraic subset, are PSC w.r.t. $\bar F$ (see Example \ref{exp-semi}).  Nevertheless, in this case, PTSC (in the complex field) can still provide some useful information for controllability preservation. For example, PTSC can guarantee that almost all real-valued realizations in ${\bf CS}(\bar A, \bar B)$ can preserve controllability under any real-valued perturbations prescribed by ${\bar F}$, noting that all such perturbations belong to ${\bf S}_{\bar F}$. On the other hand, for almost all real-valued realizations in ${\bf CS}(\bar A, \bar B)$, PSSC is a necessary condition for the existence of real-valued perturbations in ${\bf S}_{\bar F}$ that can make the perturbed systems uncontrollable. 
					
					
					\begin{example}[PTC in real field]\label{exp-semi} Consider {\footnotesize$A=\left[
							\begin{array}{cc}
								a_{11} & a_{12} \\
								a_{21} & a_{22} \\
							\end{array}
							\right]
							$, $b=\left[
							\begin{array}{c}
								0 \\
								b_2 \\
							\end{array}
							\right]
							$,} and {\footnotesize $[\Delta A, \Delta b]=\left[
							\begin{array}{ccc}
								0 & 0 & \Delta b_1 \\
								0 & 0 & 0 \\
							\end{array}
							\right]
							$}. Then, $\det{\cal C}(A+\Delta A,b+\Delta b)=a_{21}\Delta b_1^2 + (a_{22}b_2 - a_{11}b_2)\Delta b_1 - a_{12}b_2^2$. Let $g(z)=(a_{22}b_2 - a_{11}b_2)^2+4a_{21}a_{12}b_2^2$ with $z\doteq (a_{11},a_{12}, a_{21},a_{22},b_2)$. It turns out that, if $z\in \{z\in {\mathbb R}^5: g(z)\ge 0, a_{12}b_2\ne 0\}$, then there is $\Delta b_1\in {\mathbb R}$ making $(A+\Delta A,b+\Delta b)$ uncontrollable; otherwise, for $z\in \{z\in {\mathbb R}^5: g(z)< 0, a_{12}b_2\ne 0\}$, there is no $\Delta b_1\in {\mathbb R}$ making $(A+\Delta A,b+\Delta b)$ uncontrollable. 
					\end{example}

					\section{Necessary and Sufficient Condition For Single-Input Systems} \label{sec-single}
					In this section, we present a necessary and sufficient condition for PTSC in the single-input case.
					\subsection{One-edge Preservation Principle}
					At first, an one-edge preservation principle is given as follows, which is crucial to our subsequent derivations.
					
					\begin{proposition}[One-edge preservation principle] \label{OneEdgeEquivalence} Suppose $(\bar A, \bar b)$ is structurally controllable. $(\bar A,\bar b)$ is PSSC w.r.t. $\bar F$, if and only if there is one edge $e\in {\cal E}_{\bar F}$, such that $[\bar A, \bar b]\vee \bar F^{\{e\}}$ is PSSC w.r.t. ${\bar F}_{\{e\}}$, where $\bar F^{\{e\}}$ denotes the structured matrix associated with the graph ${\cal G}(\bar F)-\{e\}$, and ${\bar F}_{\{e\}}$ the structured matrix obtained from ${\bar F}$ by preserving only the entry corresponding to $e$. 
					\end{proposition}
					
					\begin{proof} Let $p\doteq (p_1,...,p_r)$ and $\bar p=(\bar p_1,...,\bar p_l)$ be such that are defined in the proof of Proposition \ref{generic_PSC_single}. From the analysis in that proof, $(\bar A, \bar b)$ is PSSC w.r.t. $\bar F$, if and only if there exists a $\bar p_j$, $j\in \{1,...,l\}$, that has a degree no less than one in $\det {\cal C}(A+\Delta A, b+\Delta b)$ (expressed in (\ref{determinant_of_Ab})). Let $e$ be the edge corresponding to $\bar p_j$. Suppose that the coefficient of $\bar p_j^k$ is nonzero for some degree $k\ge 1$. Since the coefficient of $\bar p_j^k$ in $\det {\cal C}(A+\Delta A, b+\Delta b)$ is a polynomial of $p$ and ${\bar p}\backslash \{\bar p_j\}$, it always equals the coefficient of $\bar p_j^k$ in $\det {\cal C}(A'+\Delta A', b'+\Delta b')$, where $(A',b')$ has the system graph ${\cal G}(\bar A, \bar b)\cup {\cal G}({\bar F})- \{e\}$, and $[\Delta A',\Delta b']$ corresponds to the perturbation ${\bar F}_{\{e\}}$, noting that $[A+\Delta A, b+\Delta b]\equiv[A'+\Delta A', b'+\Delta b']$ in the symbolic operation sense. Upon observing this, the proposed statement follows immediately.
					\end{proof}
					
					It is remarkable that the one-edge preservation principle does {\emph{not}} mean perturbing only one entry is enough to destroy controllability. Instead, it means we can regard $||\bar F||_0-1$ perturbed entries of $\bar F$ as unknown generic entries ({\emph{i.e., their values can be chosen randomly; but not fixed zero}}) and find suitable values for the last entry. This principle indicates that for single-input systems, verifying the PTSC w.r.t. an arbitrary perturbation structure can reduce to an equivalent problem with a single-edge perturbation structure. Having observed this, in the following, we will first give the conditions for the absence of zero uncontrollable modes and nonzero uncontrollable modes, respectively, in the single-edge perturbation scenario. Recall that an uncontrollable mode for $(A,B)$ is a $\lambda \in {\mathbb C}$ making ${\rm rank}([A-\lambda I, B])<n$. 
					
					\subsection{Conditions for Zero/Nonzero Uncontrollable Modes} \label{sec-nonzero}
					Let $\bar H\doteq [\bar A, \bar b]$. For $j\in {\cal J}_{n+1}$, let $r_j={\rm grank}(\bar H[{\cal J}_{n},{\cal J}_{n+1}\backslash \{j\}])$. Define  sets ${\cal I}_j$ and ${\cal I}^*_j$ as ${\cal I}_j=\left\{{\cal I}\subseteq {\cal J}_n: {\rm grank} (\bar H[{\cal I},{\cal J}_{n+1}\backslash \{j\}])=r_j, |{\cal I}|=r_j\right\}$, ${\cal I}^{*}_j=\left\{{\cal J}_n\backslash w: w\in {\cal I}_j\right\}. $ Based on these definitions, the following proposition gives a necessary and sufficient condition for the absence of zero uncontrollable modes in the single-edge perturbation scenario.
					
					\begin{proposition}\label{ZeroCondition} Suppose that $(\bar A, \bar b)$ is structurally controllable, and there is only one nonzero entry in ${\bar F}$ with its position being $(i,j)$. Then, for almost all $(A,b)\in {\bf CS}(\bar A,\bar b)$, there is no $[\Delta A, \Delta b]\in {\bf S}_{{\bar F}}$ such that a nonzero $n$-vector $q$ exists making $q^{\intercal}[A+\Delta A, b+\Delta b]=0$, if and only if $i\notin {\cal I}^{*}_j$.
					\end{proposition}

					To prove Proposition \ref{ZeroCondition}, we need the following lemma.
					
					\begin{lemma}\label{basiclemma} Given a matrix $H\in {\mathbb C}^{p\times q}$ of rank $p-1$, let $x\in {\mathbb C}^{p}$ be a nonzero vector in the left null space of $H$. Then, for any $i\in {\cal J}_p$,  $x_i\ne 0$, if and only if $H[{{\cal J}_p\backslash \{i\},{\cal J}_q}]$ is  of full row rank.
					\end{lemma}

					
					\begin{proof} 
						Without losing any generality, consider $i=1$. Let $h_1=H[\{1\},{\cal J}_q]$, $h_{2:p}=H[{\cal J}_p\backslash \{1\},{\cal J}_q]$, and accordingly, let $x_{2:p}=[x_2,...,x_p]$. By definition,
						\begin{equation}\label{nulleq} x^\intercal H=[x_1,x_{2:p}][h^\intercal_1,h^\intercal_{2:p}]^\intercal=x_1h_1+x_{2:p}h_{2:p}=0.\end{equation}

						{\bf If}:  Suppose $h_{2:p}$ is of full  row rank but $x_1=0$. Then, from (\ref{nulleq}), $x_{2:p}h_{2:p}=0$. As $x_2\ne 0$ (otherwise $x=0$), $h_{2:p}$ is of row rank deficient, causing a contradiction.	{\bf Only if}: Suppose $x_1\ne 0$. If $h_{2:p}$ is of row rank deficient, then there is $y\in {\mathbb C}^{p-1}$ with $y\ne0$ making $y^{\intercal} h_{2:p}=0$. Hence, $[0,y^{\intercal}]^\intercal$ is also in the left null space of $H$. As $x_1\ne 0$, this contradicts the fact that $H$ has rank $p-1$.
					\end{proof}
					

					{\bf {\emph {Proof of Proposition \ref{ZeroCondition}:}}} Sufficiency:  Since $(\bar A, \bar b)$ is structurally controllable, ${\rm grank}(\bar H)= n$, which means $r_j=n$ or $n-1$. If $r_j=n$, then ${\cal I}^*_j=\emptyset$, leading to that no $q(\ne 0)$ exists making $q^{\intercal}[A+\Delta A, b+\Delta b]=0$ for almost all $(A,b)\in {\bf CS}(\bar A,\bar b)$. Now suppose $r_j=n-1$. A vector $q (\ne 0)$ making $q^{\intercal}[A+\Delta A, b+\Delta b]=0$ must lie in the left null space of $H[{\cal J}_{n}, J_{n+1}\backslash \{j\}]$, for almost all $H\in {\bf S}_{\bar H}$. As $r_j=n-1$, for almost all $(A,b)\in {\bf CS}(\bar A,\bar b)$, ${\cal I}^*_j$ consists of all the nonzero positions of $q$ according to Lemma \ref{basiclemma}.  As a result, if $i\notin {\cal I}^*_j$, for all $[\Delta A, \Delta b]\in {\bf S}_{{\bar F}}$,
					$q^{\intercal}[A+\Delta A, b+\Delta b][{\cal J}_{n},\{j\}]=\sum \nolimits_{k\in {\cal I}^*_j}q_k[A,b]_{kj}\ne 0,$
					where the inequality is due to the fact that $[A,b]$ has full row rank.
					
					Necessity: Assume that $i\in {\cal I}^*_j$. As ${\cal I}^*_j \ne \emptyset$ and $(\bar A, \bar b)$ is structurally controllable, $[A,b][{\cal J}_n,{\cal J}_{n+1}\backslash \{j\}]$ has rank $n-1$ for all $(A,b)\in {\bf CS}(\bar A,\bar b)$. Let $q$ be a nonzero vector in the left null space of  $[A,b][{\cal J}_n,{\cal J}_{n+1}\backslash \{j\}]$.
					According to Lemma \ref{basiclemma}, as $i\in {\cal I}^*_j$, we have $q_{i}\ne 0$ for almost all $(A,b)\in {\bf CS}(\bar A,\bar b)$. By setting $[\Delta A, \Delta b]_{ij}=-1/q_i\sum \nolimits_{k\in {\cal I}^*_j}q_k[A,b]_{kj}$, we get $q^{\intercal}[A+\Delta A, b+\Delta b][{\cal J}_n,\{j\}]= q_i[\Delta A,\Delta b]_{ij}+\sum \nolimits_{k\in {\cal I}^*_j}q_k[A,b]_{kj}=0,$
					which makes $q^{\intercal}[A+\Delta A, b+\Delta b]=0$.   $\hfill \Box$
					
					\begin{remark}\label{equal-pro1} From the proof above, $i\notin {\cal I}_j^*$ is equivalent to that, ${\rm grank}(\bar H[{\cal J}_n, {\cal J}_{n+1}\backslash \{j\}])=n$ (corresponding to $r_j=n$) or ${\rm grank}(\bar H[{\cal J}_n\backslash \{i\}, {\cal J}_{n+1}\backslash \{j\}])=n-2$ (corresponding to $r_j=n-1$ but ${\rm grank}(\bar H[{\cal J}_n\backslash \{i\}, {\cal J}_{n+1}\backslash \{j\}])<n-1$). Moreover, since adding a row to a matrix can increase its rank by at most one,  the latter two conditions are mutually exclusive.
					\end{remark}
				
					Next, we present a necessary and sufficient condition for the presence of nonzero uncontrollable modes using the DM-decomposition.	
					For $j\in {\cal J}_{n+1}$, let $j_{\rm c}\doteq {\cal J}_{n+1}\backslash \{j\}$. Moreover, define a generic matrix pencil as $H_\lambda\doteq [A-\lambda I,b]$, $H_\lambda^j\doteq H_{\lambda}[{\cal J}_n,\{j\}]$, and $H_\lambda^{{j_{\rm c}}}\doteq H_{\lambda}[{\cal J}_n,j_c]$, in which $[A,b]$ is a generic realization of $[\bar A, \bar b]$. Here, the subscript $\lambda$ indicates a matrix-valued function of $\lambda$.  Let ${\cal B}(H_\lambda)=({\cal V}^+,{\cal V}^-,{\cal E})$ be the bipartite graph associated with $H_\lambda$, where ${\cal V}^+=\{x_1,...,x_n\}$, ${\cal V}^-
					=\{v_1,...,v_{n+1}\}$, and ${\cal E}=\{(x_i,v_k): {\cal E}_I\cup {\cal E}_{[\bar A,\bar b]}\}$ with ${\cal E}_I=\{(x_i,v_i):i=1,...,n\}$, ${\cal E}_{[\bar A,\bar b]}=\{(x_i,v_k):
					[\bar A,\bar b]_{ik}\ne 0\}$. No parallel edges are included even if ${\cal E}_I\cap {\cal E}_{[\bar A,\bar b]}\ne \emptyset$. An edge is called a $\lambda$-edge if it belongs to ${\cal E}_I$,  and a self-loop if it belongs to ${\cal E}_I\cap {\cal E}_{[\bar A,\bar b]}$. Note by definition, a self-loop is also a $\lambda$-edge. Let ${\cal B}(H_\lambda^{{j_{\rm c}}})$ be the bipartite graph associated with $H_\lambda^{{j_{\rm c}}}$, that is, ${\cal B}(H_\lambda^{{j_{\rm c}}})= {\cal B}(H_{\lambda})-\{v_j\}$.

					\begin{lemma} \label{fullrank}
						Suppose $(\bar A,\bar b)$ is structurally controllable. Then ${\rm mt}({\cal B}(H_\lambda^{{j_{\rm c}}}))=n$ for all $j\in {\cal J}_{n+1}$.
					\end{lemma}
					
					\begin{proof} If $j=n+1$, it is obvious ${\rm mt}({\cal B}(H_\lambda^{{j_{\rm c}}}))=n$ as ${\cal E}_I$ is a matching with size $n$. Now consider $j\in\{1,...,n\}$. As  $(\bar A,\bar b)$ is structurally controllable, from Lemma \ref{theo-strucon}, there is a path from $x_{n+1}$ to $x_j$ in the system graph ${\cal G}(\bar A, \bar b)$. Denote such a path by $\{(x_{n+1},x_{j_1}),(x_{j_1},x_{j_2}),...,(x_{j_{r-1}},x_{j_r})\}$ with  $\{j_1,...,j_r\}\subseteq {\cal J}_n$ and $j_r=j$. Since each $(x_{j_k},x_{j_{k+1}})$ in ${\cal G}(\bar A, \bar b)$ corresponds to $(x_{j_{k+1}},v_{j_k})$ in ${\cal B}(H_\lambda)$, $\{(x_{j_1},v_{n+1}),(x_{j_2},v_{j_1}),...,(x_{j_r},v_{j_{r-1}})\}\cup \{(x_i,v_i): i\in {\cal J}_n\backslash \{j_1,...,j_r\}\}$ forms a matching with size $n$ in ${\cal B}(H_\lambda^{{j_{\rm c}}})$.
					\end{proof}
					
					Let ${\cal G}^{j_{\rm c}}_k=({\cal V}^+_k,{\cal V}_k^-,{\cal E}_k)$ ($k=0,1,...,d,\infty$) be the DM-components of ${\cal B}(H_\lambda^{{j_{\rm c}}})$. From Lemma \ref{fullrank}, we know that both the horizontal tail and the vertical one are empty. Let $M^{j_{\rm c}}_\lambda$ be the DM-decomposition of $H_\lambda^{{j_{\rm c}}}$ with  permutation matrices $P$ and $Q$, i.e.,
					\begin{equation}\label{DMdecomposition} PH_\lambda^{{j_{\rm c}}}Q=\left[
						\begin{array}{ccc}
							M^{j_{\rm c}}_1(\lambda) & \cdots & M^{j_{\rm c}}_{1d}(\lambda) \\
							0 & \ddots & \vdots \\
							0 & \cdots & M^{j_{\rm c}}_d(\lambda)  \\
						\end{array}
						\right]\doteq M^{j_{\rm c}}_\lambda.\end{equation}
					Moreover, define $M_\lambda^j\doteq PH_\lambda^j$. Suppose that $x_i$ is the $\bar i$th vertex in ${\cal V}^+_{i^*}$ ($1\le \bar i\le |{\cal V}^+_{i^*}|$, $1\le i^* \le d$).
					For $k\in \{1,...,d\}$, let $\gamma_{\min}({\cal G}_k^{j_{\rm c}})$ and $\gamma_{\max}({\cal G}_k^{j_{\rm c}})$ be respectively the minimum number of $\lambda$-edges and maximum number of $\lambda$-edges contained in a matching over all maximum matchings of ${\cal G}_k^{j_{\rm c}}$. 
					 Afterwards, define a boolean function $\gamma_{\rm nz}(\cdot)$ for ${\cal G}_k^{j_{\rm c}}$ as
					\begin{equation}\label{fun-nz} \gamma_{\rm nz}({\cal G}_k^{j_{\rm c}})=
						\begin{cases} 1 & {\begin{array}{c} \text{if}\ \gamma_{\max}({\cal G}_k^{j_{\rm c}})-\gamma_{\min}({\cal G}_k^{j_{\rm c}})>0 \\ \text{or} \ {\cal G}_k^{j_{\rm c}} \ {\text{contains a self-loop}} \end{array}} \\
							0 &  \text{otherwise}.
					\end{cases}\end{equation}By assigning the weight $W_k: {\cal E}_k \rightarrow \{0,1\}$ for ${\cal G}_k^{j_{\rm c}}$ as: $W_k(e)=1$ if $e\in {\cal E}_k$ is a $\lambda$-edge, and $W_k(e)=0$ otherwise.
It is then obvious that $\gamma_{\max}({\cal G}_k^{j_{\rm c}})$ equals ${\rm maxWMM}({\cal G}_k^{j_{\rm c}})$, and $\gamma_{\min}({\cal G}_k^{j_{\rm c}}) = {\rm minWMM({\cal G})}$. Hence, $\gamma_{\rm nz}({\cal G}_k^{j_{\rm c}})$ can be determined in polynomial time via the maximum weighted matching algorithms \cite{Murota_Book}.
					From Lemma \ref{matrix_pencil} in the appendix, $\gamma_{\rm nz}({\cal G}_k^{j_{\rm c}})=1$ means $\det (M_k^{j_c}(\lambda))$ generically has nonzero roots for $\lambda$, while $\gamma_{\rm nz}({\cal G}_k^{j_{\rm c}})=0$ means the contrary, $k\in \{1,...,d\}$.

%
					
					
					Furthermore, define a set
					\begin{equation}\label{important-set} {\Omega}_j=\{k\in {\mathbb N}: 1\le k\le i^*, {\gamma}_{\rm nz}({\cal G}^{j_{\rm c}}_k)= 1\}.\end{equation}
					For each $k\in {\Omega}_j$, define a weighted bipartite graph ${\cal G}^{j_{\rm c}}_{ki^*}=\left(\bar {\cal V}^+_{ki^*}, \bar {\cal V}^-_{ki^*},  \bar {\cal E}_{ki^*}, W \right)$, where $\bar {\cal V}^+_{ki^*}={\cal V}^{+}_k\cup {\cal V}^{+}_{k+1}\cdots\cup {\cal V}^+_{i^*}\backslash \{x_i\}$, $\bar {\cal V}^-_{ki^*}= {\cal V}^{-}_k\cup {\cal V}^{-}_{k+1}\cdots \cup {\cal V}^-_{i^*} $, $\bar {\cal E}_{ki^*}= \{(x_q,v_l)\in {\cal E}: x_q\in \bar {\cal V}^+_{ki^*}, v_l\in \bar {\cal V}^-_{ki^*}\}$, and the weight $W(e): \bar {\cal E}_{ki^*} \rightarrow \{0,1\}$
					$$W(e)=
					\begin{cases} 1 & \text{if}\ e\in {\cal E}_{k}\\
						0 & {\text{otherwise.}}
					\end{cases}$$
					
					\begin{proposition} \label{mainpro}
						\label{nonzeroCondition} Suppose $(\bar A, \bar b)$ is structurally controllable, and there is only one nonzero entry in ${\bar F}$ with its position being $(i,j)$.  Then, for almost all $(A,b)\in {\bf CS}(\bar A,\bar b)$, there is a $[\Delta A, \Delta b]\in {\bf S}_{{\bar F}}$ such that a nonzero $n$-vector $q$ exists making $q^{\intercal}[A+\Delta A-\lambda I, b+\Delta b]=0$ for some nonzero $\lambda\in {\mathbb C}$, if and only if there exists a $k\in {\Omega}_j$ such that ${\rm minWMM}({\cal G}^{j_{\rm c}}_{ki^*})< |{\cal V}^{+}_k|$, with ${\cal G}^{j_{\rm c}}_{ki^*}$ defined above.
						
					\end{proposition}

					The proof relies on a series of nontrivial results on the roots of determinants of generic matrix pencils, which is postponed to Appendix~\ref{mainpro-proof}.
					\subsection{Necessary and Sufficient Condition}
					We are now giving a necessary and sufficient condition for PTSC with general perturbation structures.
					
					\begin{theorem}\label{NeceSuf} Consider a structurally controllable pair $(\bar A,\bar b)$ and the perturbation structure $\bar F$. For each edge $e\doteq (x_j,x_i)\in {\cal E}_{\bar F}$, let $[\bar A^e,\bar b^e]=[\bar A, \bar b]\vee \bar F^{\{e\}}$, with $\bar F^{\{e\}}$ defined in Proposition \ref{OneEdgeEquivalence}. Moreover, let ${\Omega}_j$ and ${\cal G}^{j_{\rm c}}_{ki^*}$ be defined in the same way as in Proposition \ref{nonzeroCondition}, in which $(\bar A, \bar b)$ shall be replaced with $(\bar A^e,\bar b^e)$. Then, $(\bar A,\bar b)$ is PTSC w.r.t. $\bar F$, if and only if for each edge $e\doteq (x_j,x_i)\in {\cal E}_{\bar F}$, it holds simultaneously:
						
						1) ${\rm grank}(\bar H[{\cal J}_n, {\cal J}_{n+1}\backslash \{j\}])=n$ or ${\rm grank}(\bar H[{\cal J}_n\backslash \{i\}, {\cal J}_{n+1}\backslash \{j\}])=n-2$, with $\bar H=[\bar A^e,\bar b^e]$;
						
						2) ${\Omega}_j=\emptyset$, or for each $k\in {\Omega}_j$, ${\rm minWMM}({\cal G}^{j_{\rm c}}_{ki^*})= |{\cal V}^{+}_k|$.
					\end{theorem}
					
					
					\begin{proof}The result is immediate from the PBH test and Propositions \ref{OneEdgeEquivalence}-\ref{nonzeroCondition}.
					\end{proof}
					
					Since each step in Theorem \ref{NeceSuf} can be implemented in polynomial time, its verification has polynomial time complexity, too. To be specific, for each edge $e\in {\cal E}_{\bar F}$, to verify condition 1), we can invoke the Hopcroft-Karp algorithm \cite{Geor1993Graph} twice, which incurs time complexity $O(n^{0.5}|{\cal E}_{[\bar A, \bar b]}\cup {\cal E}_{\bar F}|)$ $\to O(n^{2.5})$. As for condition 2), the DM-decomposition incurs $ O(n^{2.5})$, and computing ${\rm minWMM}({\cal G}^{j_{\rm c}}_{ki^*})$ costs $O(n^3)$ \cite{DB_West_graph}. Since $|\Omega_j|\le n$, for each $e\in {\cal E}_{\bar F}$, verifying condition 2) takes at most $O(n^{2.5}+n*n^3)$. To sum up, verifying Theorem \ref{NeceSuf} incurs time complexity at most $O(|{\cal E}_{\bar F}|(n^{2.5}+n^4))$, i.e., $O(|{\cal E}_{\bar F}|n^4)$. 
					
					
					
					\begin{example}[Example \ref{exp1} contd]\label{exp1-con} Let us revisit Example \ref{exp1}. Consider the perturbation $[\Delta A_2, \Delta b_2]$. For edge $e=(x_5,x_4)$,  the DM-decomposition of $H_{\lambda}^{j_{\rm c}}$ ($j=5$) associated with $(\bar A^e,\bar b^e)$ and the corresponding $M_{\lambda}^j$ are respectively
						$$ M_{\lambda}^{j_{\rm c}}\!=\!\left[
						\begin{array}{c|c|c|c}
							\Delta a_{33}-\lambda  &    0     &     a_{32}   &    0     \\
							\hline
							& a_{44}-\lambda  &    a_{42}    &    a_{41}       \\
							\cline{2-4} \multicolumn{2}{c|}{} & -\lambda & a_{21} \\
							\cline{3-4}        \multicolumn{3}{c|}{} &  -\lambda
						\end{array}
						\right], M_{\lambda}^{j}\!=\!\left[\!
						\begin{array}{c}
							0 \\
							\Delta b_4 \\
							0 \\
							b_1 \\
						\end{array}
						\!\right].
						$$
						It can be obtained that, $i^*=2$, and $\Omega_j=\{1,2\}$. Since $i^*\in \Omega_j$, according to Proposition \ref{nonzeroCondition}, the corresponding perturbed system can have nonzero uncontrollable modes (in fact, if $i^*\in \Omega_j$, then the condition in Proposition \ref{nonzeroCondition} is automatically satisfied). Therefore, $(\bar A,\bar b)$ is PSSC w.r.t. $[\Delta \bar A_2,\Delta \bar b_2]$, which is consistent with Example \ref{exp1}. On the other hand, consider the perturbation $(\Delta A_1, \Delta b_1)$. For the edge $e=(x_4,x_1)$, upon letting $j=4$, we obtain ${\cal I}_j^*=\{2,3,4\}$, which means $1\notin {\cal I}_j^*$. Hence, the condition in  Proposition \ref{ZeroCondition} is satisfied. Moreover, the associated $M_{\lambda}^{j_{\rm c}}$ and $M_{\lambda}^{j}$ are respectively
						$$M_{\lambda}^{j_{\rm c}}=\left[
						\begin{array}{c|c|c|c}
							b_1  &    \Delta a_{13}     &     -\lambda   &    0     \\
							\hline
							& -\lambda  &   0    &    a_{32}       \\
							\cline{2-4} \multicolumn{2}{c|}{} & a_{41} & a_{42} \\
							\cline{3-4} \multicolumn{2}{c|}{} & a_{21} & -\lambda
						\end{array}
						\right], M_{\lambda}^{j}=\left[
						\begin{array}{c}
							\Delta a_{14} \\
							0 \\
							a_{44}-\lambda \\
							0 \\
						\end{array}
						\right], $$from which, $i^*=1$, and $\Omega_j=\emptyset$. It means condition 2) of Theorem \ref{NeceSuf} is satisfied.
						Similar analysis could be applied to the edge $e=(x_3,x_1)$, and it turns out that both conditions in Theorem \ref{NeceSuf} hold. Therefore, $(\bar A,\bar b)$ is PTSC w.r.t. $[\Delta \bar A_1,\Delta \bar b_1]$, which is also consistent with Example \ref{exp1}. 
					\end{example}
					
				{	\begin{remark}
						The techniques in this section could also lead to some necessary conditions for PTSC of multi-input systems. However, due to that the one-edge preservation principle does not hold for general multi-input systems, complete criteria for the multi-input case need further investigation. 
					\end{remark} }
					
					\section{Graph-Theoretic Criteria for Single-Input Case} \label{graph-theoretic}
					In this section, we give some intuitive graph-theoretic conditions for PSSC in the single-input case. In particular, we will discuss PSSC w.r.t. edge addition, relate it to the {\emph{cactus}}, and present some graphical criteria.  
					
					
					
					\begin{proposition}[Edge addition]\label{property-unit}  Consider two single-input systems $(\bar A_1,\bar b_1)$, $(\bar A_2,\bar b_2)$ and a perturbation structure ${\bar F}\in \{0,*\}^{n\times (n+1)}$. Suppose that $[\bar A_1,\bar b_1]\subseteq [\bar A_2,\bar b_2]$. If $(\bar A_1,\bar b_1)$ is PSSC w.r.t. $\bar F$, then $(\bar A_2,\bar b_2)$ is PSSC w.r.t. $\bar F$.
					\end{proposition}
					
					\begin{proof} Suppose $(\bar A_1,\bar b_1)$ is PSSC w.r.t. $\bar F$. By Proposition \ref{OneEdgeEquivalence}, there is one edge $e\in {\cal E}_{\bar F}$, such that $[\bar A_1,\bar b_1]\vee \bar F^{\{e\}}$ is PSSC w.r.t. ${\bar F}_{\{e\}}$. Let $p_1,...,p_t$ be values for the nonzero entries in a generic realization $[A'_2,b'_2]$ of $[\bar A_2, \bar b_2]\cup \bar F^{\{e\}}$, and $p\doteq (p_1,...,p_t)$. Moreover,  let $s$ be the value of the unique nonzero entry in $[\Delta A', \Delta b']\in {\bf S}_{{\bar F}_{\{e\}}}$. Rewrite $\det{\cal C}(A'_2+\Delta A', b'_2+\Delta b')$ as a polynomial of $s$ as
						$\det{\cal C}(A'_2+\Delta A', b'_2+\Delta b')=f_r(p)s^r+\cdots+f_1(p)s+f_0(p),$
						where $f_r(p),\cdots, f_0(p)$ are polynomials of $p$, $r\ge 0$. As $[\bar A_1,\bar b_1]\vee \bar F^{\{e\}}\subseteq [\bar A_2,\bar b_2]\vee \bar F^{\{e\}}$, $(A'_2, b'_2)$ is generically controllable, which means $f_0(p)\not\equiv 0$ (by letting $s=0$).
						
						Now let $\bar p$ be obtained from $p$ by making the $p_i$ elements of the nonzero entries in $[\bar A_2,\bar b_2]$ but not in $[\bar A_1,\bar b_1]$  be zero and preserving the rest. For the generic realization $[A'_1,b'_1]$ of $[\bar A_1, \bar b_1]\vee \bar F^{\{e\}}$ with parameters $\bar p$, we have $\det{\cal C}(A'_1+\Delta A', b'_1+\Delta b')=f_r(\bar p)s^r+\cdots+f_0(\bar p)$. As $[\bar A_1,\bar b_1]\vee \bar F^{\{e\}}$ is PSSC w.r.t. ${\bar F}_{\{e\}}$, from the proof of Proposition \ref{generic_PSC_single}, there must exist one $t\in \{1,...,r\}$ such that $f_t(\bar p)\not\equiv 0$.
						This requires $f_t(p)\not\equiv 0$. Thereby, $[\bar A_2,\bar b_2]\vee \bar F^{\{e\}}$ is PSSC w.r.t. $\bar F_{\{e\}}$. By the one-edge principle principle, $(\bar A_2, \bar b_2)$ is PSSC w.r.t. $\bar F$. 
					\end{proof}
					
					Proposition \ref{property-unit} simply means adding extra unknown generic entries to a single-input system which is PSSC will not make the resulting system PTSC (w.r.t. the same $\bar F$). Recall that a $\emph{cactus}$ is the minimal structure that can preserve structural controllability (in single-input case), in the sense that removing any edge from this structure will result in a structurally uncontrollable structure \cite{liu2016control}. Based on these observations, the following result gives an intuitive sufficient condition for PSSC.
					
					\begin{proposition} \label{cactus-pro} Consider a structurally controllable pair $(\bar A, \bar b)$ and $\bar F$. If there exists an edge $e\in {\cal E}_{\bar F}$ that is contained in a cactus of ${\cal G}(\bar A,\bar b)\cup {\cal G}(\bar F)$, then $(\bar A, \bar b)$ is PSSC w.r.t. $\bar F$.
					\end{proposition}
					
					\begin{proof}
						 Let $[\bar A_1,\bar b_1]\subseteq [\bar A, \bar b]$ and $\bar F_1 \subseteq \bar F$ be such that ${\cal G}(\bar A_1,\bar b_1)\cup {\cal G}(\bar F_1)$ corresponds exactly to the cactus of ${\cal G}(\bar A,\bar b)\cup {\cal G}(\bar F)$ that contains $e$. Since a cactus is the minimal structure for structural controllability, upon letting the entry of $[\bar A_1,\bar b_1]\vee \bar F_1$ corresponding to $e$ be fixed zero, the resulting system will be structurally uncontrollable. This means, $[\bar A_1,\bar b_1]\vee \bar F_1$ is PSSC w.r.t. $\bar F_{\{e\}}$, with $\bar F_{\{e\}}$ defined in Proposition \ref{OneEdgeEquivalence}. By Proposition \ref{property-unit}, $[\bar A, \bar b] \vee \bar F$ is PSSC w.r.t. $\bar F_{\{e\}}$. As $(\bar A, \bar b)$ is structurally controllable, by Definition \ref{DefinitionForPTSC}, $(\bar A, \bar b)$ is PSSC w.r.t. $\bar F$.
					\end{proof}

					Proposition \ref{cactus-pro} provides only a sufficient condition for PSSC. The perturbed edge $e$ leading to PSSC is not necessarily contained in a cactus, as seen from the following example.
					
					\begin{example}[Cactus is not necessary] \label{cactus-not-necessary}
						Consider $[A,b]\in {\bf S}_{[\bar A, \bar b]}$ as
						$${\footnotesize 
							A=\left[
						\begin{array}{cccc}
							0 &     0 &     0 &     0 \\
							a_{21}  &   0 &     0 &     0 \\
							0 &     a_{32}  &   0 &   0 \\
							a_{41} &    0 &     0 &     a_{44}
						\end{array}
						\right],b=\left[
						\begin{array}{c}
							b_1 \\
							0 \\
							0 \\
							0 \\
						\end{array}
						\right]}.
						$$
						${\bar F}$ has only one nonzero entry in its $(3,3)$th position, being $\Delta a_{33}$. The system graph ${\cal G}(\bar A,\bar b)\cup {\cal G}(\bar F)$, as well as its cactus, is given in Fig. \ref{cactus-1}. It is seen that edge $(x_3,x_3)$ is not contained in the cactus. However, through some simple calculations,
						$\det{\cal C}(A+\Delta A,b+\Delta b)=a_{21}^2a_{32}a_{41}b_1^4(a_{44}^2- \Delta a_{33}a_{44})$. Letting $\Delta a_{33}=a_{44}$, $(A+\Delta A, b+\Delta b)$ becomes uncontrollable. This means, $(\bar A, \bar b)$ is PSSC w.r.t. $\bar F$.
					\end{example}
					
					
					\begin{figure}
						\centering
						\includegraphics[width=1.9in]{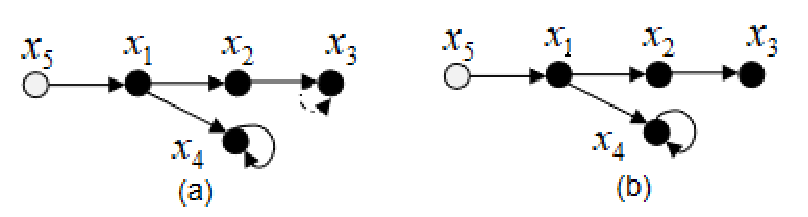}\\
						\caption{${\cal G}(\bar A,\bar b)\cup {\cal G}(\bar F)$ (a) and its cactus (b) in Example \ref{cactus-not-necessary}. Dotted edges represent perturbed edges. }\label{cactus-1}
					\end{figure}
					
					In the following, we provide separately equivalent graphical forms the two conditions in Theorem \ref{NeceSuf}, which may be more intuitive to verify. Owing to the one-edge preservation principle, it suffices to consider the case where the perturbation structure $\bar F$ contains only one nonzero entry, and suppose its position is $(i,j)$. The case where $\bar F$ contains more nonzero entries than one can be obtained directly (in the same way as Theorem \ref{NeceSuf}).
					
					
					To this end, some terminologies are first introduced from \cite{generic}. A number of vertex-disjoint elementary paths (resp. cycles) is called a {\emph{path family}} (resp. {\emph{cycle family}}). The union of vertex-disjoint path families and cycle families is called a {\emph{path-cycle family}}. The size of a path-cycle family is the number of edges it contains.
					Next, construct an auxiliary graph associated with $(\bar A, \bar b)$ as $\tilde {\cal G}(\bar A, \bar b)=(\tilde {\cal X} \cup \tilde {\cal U}, \tilde {{\cal E}}_{\bar A} \cup \tilde {{\cal E}}_{\bar b} \cup \tilde {{\cal E}}_{I})$, with $\tilde {\cal X}=\{x_1^{\rm i}, x_1^{\rm o},\cdots, x_n^{\rm i}, x_n^{\rm o}\}$, $\tilde {\cal U}=\{x_{n+1}^{\rm i}, x_{n+1}^{\rm o}\}$, $\tilde {\cal E}_{\bar A}=\{(x_k^{\rm o}, x_l^{\rm i}): \bar A_{lk}\ne 0\}$, $\tilde {\cal E}_{\bar b}=\{(x_{n+1}^{\rm o}, x_k^{\rm i}): \bar b_k\ne 0\}$, and $\tilde {\cal E}_{I}=\{(x_k^{\rm o},x_k^{\rm i}): k=1,...,n\}$. Simply speaking, $\tilde {\cal G}(\bar A, \bar b)$ is obtained from ${\cal G}(\bar A, \bar b)\cup ({\cal X}, {\cal E}_I)$ by duplicating its vertices with a pair $\{x_k^{\rm i}, x_k^{\rm o}\}$, in which $x_k^{\rm i}$ receives the ingoing edges to $x_k$, and $x_k^{\rm o}$ sends the outgoing edges from $x_k$ ($k=1,...,n+1$). No parallel edges are included even if $\tilde {\cal {\cal E}}_{I}\cap \tilde {\cal E}_{\bar A}\ne \emptyset$. Similar to Section \ref{sec-nonzero}, a $\lambda$-edge is an edge in $\tilde {\cal E}_I$, and a self-loop is an edge in $\tilde {\cal E}_{I}\cap \tilde {\cal E}_{\bar A}$.
					Let $\cal M$ be a path family of $\tilde {\cal G}(\bar A,\bar b)- \{x_{n+1}^{\rm i},x_{j}^{\rm o}\}$ with size $n$, recalling $(i,j)$ for the only nonzero entry of $\bar F$. From Lemma \ref{fullrank}, such $\cal M$ always exists. Let $\tilde {\cal G}_{\cal M}$ be the graph obtained from $\tilde {\cal G}(\bar A,\bar b)-\{x_{n+1}^{\rm i},x_{j}^{\rm o}\}$ by adding edges $\left\{(x_l^{\rm i}, x_k^{\rm o}): (x_k^{\rm o}, x_l^{\rm i})\in {\cal M}\right\}$ to it (see Fig. \ref{example-graph} for illustration). Now let $\tilde {\cal V}_k$ ($k=1,...,d$) be the SCCs of $\tilde {\cal G}_{\cal M}$, and $\tilde {\cal G}_k=(\tilde {\cal V}_k, \tilde {\cal E}_k)$ be the subgraph of $\tilde {\cal G}(\bar A,\bar b)$ induced by $\tilde {\cal V}_k$. Let $y_{\rm nz}(\cdot)$ be defined for $\tilde {\cal G}_k$ in the similar way to (\ref{fun-nz}), in which $\tilde {\cal G}_k$ is rewritten as a bipartite graph $(\tilde {\cal V}^{+}_k,\tilde {\cal V}^-_k, \tilde {\cal E}_k)$ with $\tilde {\cal V}^{+}_k=\{x_l^{\rm o}: x_l^{\rm o}\in \tilde {\cal V}_k\}$ and $\tilde {\cal V}^{-}_k=\{x_l^{\rm i}: x_l^{\rm i}\in \tilde {\cal V}_k\}$. Define
					$$\tilde \Omega_j=\{k\in {\mathbb N}: 1\le k \le d, y_{\rm nz}(\tilde {\cal G}_k)=1\}. $$
					
					Based on the above constructions, we have the following equivalent interpretations of Propositions \ref{ZeroCondition} and \ref{mainpro}, which, along with Proposition \ref{OneEdgeEquivalence}, together give a graphical necessary and sufficient condition for a single-input system to be PSSC.

%
%
%
					
					\begin{proposition}\label{graph-zero} Under the same setting therein,  the condition in Proposition \ref{ZeroCondition} is not satisfied, if and only if the maximum size over all path-cycle families in ${\cal G}(\bar A, \bar b)-{\cal E}_{x_j}^{\rm out}$ and in ${\cal G}(\bar A, \bar b)-{\cal E}_{x_j}^{\rm out}\cup {\cal E}_{x_i}^{\rm in}$, both equals $n-1$.
					\end{proposition}
					
					\begin{proof}  Since ${\rm grank}(\bar H)=n$, one has ${\rm grank}(\bar H[{\cal J}_n, {\cal J}_{n+1}\backslash \{j\}])\ge n-1$, and ${\rm grank}(\bar H[{\cal J}_n\backslash \{i\}, {\cal J}_{n+1}\backslash \{j\}])\ge n-2$, recalling $\bar H=[\bar A, \bar b]$. Combined with Remark \ref{equal-pro1}, the condition in Proposition \ref{ZeroCondition} is not satisfied, if and only if ${\rm grank}(\bar H[{\cal J}_n, {\cal J}_{n+1}\backslash \{j\}])= n-1$, and ${\rm grank}(\bar H[{\cal J}_n\backslash \{i\}, {\cal J}_{n+1}\backslash \{j\}])= n-1$. Note that there is a one-to-one correspondence between a path-cycle family in ${\cal G}(\bar A, \bar b)-{\cal E}_{x_j}^{\rm out}$ and a matching in the bipartite graph ${\cal B}(\bar H[{\cal J}_n,{\cal J}_{n+1}\backslash \{j\})$. Hence, the maximum size of path-cycle families in ${\cal G}(\bar A, \bar b)-{\cal E}_{x_j}^{\rm out}$, is equal to the size of a maximum matching in ${\cal B}(\bar H[{\cal J}_n,{\cal J}_{n+1}\backslash \{j\})$, which equals ${\rm grank}(\bar H[{\cal J}_n,{\cal J}_{n+1}\backslash \{j\}])$, and thus should be $n-1$. Similarly, it follows that  ${\rm grank}(\bar H({\cal J}_n\backslash \{i\},{\cal J}_{n+1}\backslash \{j\}))=n-1$, if and only if the maximum matching of ${\cal G}(\bar A, \bar b)-{\cal E}_{x_j}^{\rm out}\cup {\cal E}_{x_i}^{\rm in}$ has the size $n-1$. 
					\end{proof}
					

					\begin{proposition}\label{graph-nonzero} Under the same setting therein, the condition in Proposition \ref{mainpro} is satisfied, if and only if there is a path family $\tilde {\cal M}$ of $\tilde {\cal G}(\bar A,\bar b)- \{x_j^{\rm o},x_i^{\rm i}\}$ with size $n-1$ and a $k\in \tilde \Omega_j$, such that  $\tilde {\cal M}$ contains at most $|\tilde {\cal V}_k^+|-1$ edges of ${\tilde {\cal G}}_k$.
					\end{proposition}
					
				
				\begin{proof} From the proof therein, the condition in Proposition \ref{mainpro} is equivalent to that, there is a matching in ${\cal B}(H^{j_{\rm c}}_\lambda)-\{x_i\}$ with size $n-1$ that contains $|{\cal V}_k^+|-1$ edges of ${\cal G}_k^{j_{\rm c}}$ for some $k$ satisfying ${\gamma}_{\rm nz}({\cal G}_{k}^{j_{\rm c}})=1$ (referred to as condition (d)). Indeed, if such a matching of ${\cal G}_{ki^*}^{j_c}$ satisfying the condition in Proposition \ref{mainpro} exists, then this matching together with each maximum matching of ${\cal G}^{j_c}_{t}$ for $t\in \{1,...,d\}\backslash \{k,...,i^*\}$ will form a maximum matching of ${\cal B}(H^{j_{\rm c}}_\lambda)-\{x_i\}$  satisfying condition (d). On the other hand, suppose such a matching ${\cal M}_0$ satisfying condition (d) exists. Then, from 3) of Definition \ref{DM-def}, we have $\{v\in {\cal V}^-: (x,v)\in {\cal M}_0, x\in \bar {\cal V}^+_{ki^*}\}\subseteq \bar {\cal V}^-_{ki^*}$, as otherwise there is at least one $v\in {\cal V}^-$ that is incident to two vertices in ${\cal M}_0$ (one from $\bar {\cal V}^+_{ki^*}$, and the other from $ \bigcup \nolimits_{l=i^*+1}^{d}{\cal V}^+_{l}$). Hence, a matching satisfying the condition in Proposition \ref{mainpro} is contained in ${\cal M}_0$.
					
				Observe that $\tilde {\cal G}(\bar A, \bar b)$ is bipartite with bipartitions  $\{\tilde x_1^{\rm i},...,\tilde x_{n}^{\rm i}\}$ and $\{\tilde x_1^{\rm o},...,\tilde x_{n+1}^{\rm o}\}$, and therefore, there is of one-to-one correspondence to ${\cal B}(H_{\lambda})$. In addition, $\cal M$ must be a path family of $\tilde {\cal G}(\bar A, \bar b)-\{x_{n+1}^{\rm i},x_j^{\rm o}\}$, as there is no ingoing edge to $x^{\rm i}_{n+1}$. As $|{\cal M}|=n$, there is no tail in the DM-components of $\tilde {\cal G}(\bar A, \bar b)-\{x_{n+1}^{\rm i},x_j^{\rm o}\}$. From the DM-decomposition algorithm in Page 61 of \citep[Sec 2.2.3]{Murota_Book}, each SCC $\tilde {\cal V}_k$ ($k=1,...,d$) of $\tilde {\cal G}_{\cal M}$ is actually the vertex set of a DM-component of $\tilde {\cal G}(\bar A, \bar b)-\{x_{n+1}^{\rm i},x_j^{\rm o}\}$. Hence, there is a one-to-one correspondence between $\{\tilde {\cal G}_k|_{k=1}^d\}$ and $\{{\cal G}^{j_{\rm c}}_k|_{k=1}^d\}$. The equivalence between condition (d) and the one in Proposition \ref{graph-nonzero} then follows from the facts identified above. \end{proof}  

		\begin{example}\label{example6} Consider $(\bar A,\bar b)$ and $\bar F$ in Example \ref{cactus-not-necessary} ($i=j=3$). Associated with this system, ${\cal G}(\bar A, \bar b)$, $\tilde {\cal G}(\bar A,\bar b)$, ${\tilde {\cal G}}_{\cal M}$, and $\tilde {\cal G}(\bar A, \bar b)-\{x_3^{\rm o},x_3^{\rm i}\}$  are respectively given in Figs. \ref{example-graph}(a), \ref{example-graph}(b), \ref{example-graph}(c) and \ref{example-graph}(d). From Fig. \ref{example-graph}(a), there is a path-cycle family with size $4$ in ${\cal G}(\bar A, \bar b)-{\cal E}_{x_j}^{\rm out}$. Hence, the condition in Proposition \ref{graph-zero} is not satisfied. From Fig. \ref{example-graph}(c), $\tilde {\cal G}_{\cal M}$ has four SCCs, which are $\tilde {\cal V}_1=\{x_{5}^{\rm o},x_{1}^{\rm i}\}$, $\tilde {\cal V}_2=\{x_1^{\rm o},x_{2}^{\rm i}\}$, $\tilde {\cal V}_3=\{x_2^{\rm o}, x_3^{\rm i}\}$, and $\tilde {\cal V}_4=\{x_4^{\rm o}, x_4^{\rm i}\}$. Among them, only $\tilde {\cal V}_4$ corresponds to a bipartite graph associated with which the function $y_{nz}(\cdot)=1$. Hence, $\tilde \Omega_j=\{4\}$. From Fig. \ref{example-graph}(d), there is a path family with size $3$ in $\tilde {\cal G}(\bar A,\bar b)- \{x_3^{\rm o},x_3^{\rm i}\}$ that does not contain edges of $\tilde {\cal G}_4$. Therefore, the condition of Proposition \ref{graph-nonzero} is satisfied, which indicates $(\bar A, \bar b)$ is PSSC w.r.t. $\bar F$.
		\end{example}

		\begin{figure}
		\centering
	\includegraphics[width=2.5in]{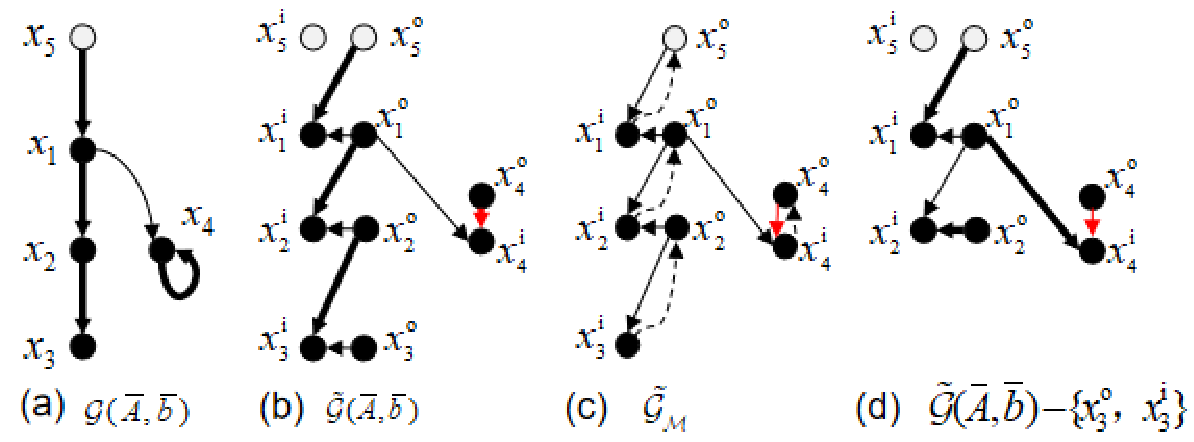}\\
		\caption{The graphs for Example \ref{example6} ($i=j=3$). In (a), bold edges represent a path-cycle family; In (b), bold edges represent a path family, while the red edge stands for a self-loop; In (c), the dotted edges are reverse edges of the path family in (b); In (d), bold edges represent a path family.}\label{example-graph}
\end{figure}

		\section{Application to SCRP} \label{apply-scrp}
		In this section, we show the application of the proposed PTSC to the SCRP.  As mentioned in Section I, the SCRP is to search the smallest perturbations (in terms of matrix norms) with a prescribed structure that result in an uncontrollable system \cite{KarowStructured2009,khare2012computing,bianchin2016observability,johnson2018structured}. Given $(A,B)$ in (\ref{plant}) and a perturbation structure $\bar F$ with compatible dimensions, one typical variant of SCRP can be formulated as follows \cite{bianchin2016observability}
		\begin{equation} \label{SCR}
			\begin{array}{l}
				\min \nolimits_{[\Delta A,\Delta B]\in {\bf S}_{\bar F}} ||[\Delta A, \Delta B]||_* \\
				{\rm s.t.}\ (A+\Delta A, B+\Delta B) \  {\rm uncontrollable}
			\end{array}
		\end{equation}in which $||\cdot||_*$ can be the Frobenius norm or $2$-norm. In some literature, an additional constraint is that $[\Delta A,\Delta B]$ needs to take real values.

		It is readily seen that PSC of $(A, B)$ w.r.t. $\bar F$ is exactly the feasibility condition for Problem (\ref{SCR}). As this property is generic for $(A,B)\in {\bf CS}(\bar A, \bar B)$, {\emph{PSSC can provide some feasibility condition for Problem (\ref{SCR}) from a generic view}}. To be specific, before looking at the exact parameters of the original system and implementing any numerical algorithms on the associated SCRP, we can check whether the corresponding structured system is PTSC w.r.t. the perturbation structure. If the answer is yes, then with probability $1$, there exist no numerical perturbations with the prescribed structure for which the corresponding perturbed system is uncontrollable; otherwise, with probability $1$, such a structured (complex-valued) numerical perturbation exists. In particular, for the single-input system, if the corresponding structured system is PTSC w.r.t. the perturbation structure, provided the original system is controllable, the associated SCRP is immediately infeasible. Furthermore, when $(\Delta A, \Delta B)$ is restricted to be real-valued, PSSC is a necessary condition for feasibility of the associated SCRP for almost all real-valued realizations in ${\bf CS}(\bar A, \bar B)$ (see Section \ref{sec_gener}).
		
		To demonstrate the above assertions, we test the performance of Algorithm 3.5 in \cite{khare2012computing}, which is an iterative algorithm for SCRP based on structured total least squares. Below is a randomly-generated realization of the structured system in Example \ref{exp1}
		{\begin{equation}\label{exampsys}{\footnotesize A\!=\!\left[\!
					\begin{array}{cccc}
						0 & 0 & 0 & 0 \\
						5.3165  &       0     &    0     &    0 \\
						0  &  9.0428    &     0      &   0 \\
						0.5454 &   6.3018  & 0 &           9.6296
					\end{array}
					\!\right],b\!=\!\left[\!
					\begin{array}{c}
						5.3401 \\
						0 \\
						0 \\
						0 \\
					\end{array}
					\!\right].}\end{equation}
			For two perturbation structures $[ \Delta \bar A_i,\Delta \bar b_i]$ ($i=1,2$) in Example \ref{exp1}, we record the corresponding  performance of Algorithm 3.5 in Table \ref{alg-perfor}, in which `complex field' (resp. `real field') means the corresponding perturbations are complex-valued (resp. real-valued). The algorithm parameters are selected as suggested in \cite{khare2012computing} (in particular, the penalty factor $w=10^{9}$ in Equation (22) of \cite{khare2012computing}, the tolerance $\varepsilon=10^{-10}$, and the maximum number of iterations is $150000$; see \cite{khare2012computing} for details). In Table \ref{alg-perfor}, `terminated tolerance' refers to the actual difference between two successive iterative parameters (corresponding to the tolerance $\varepsilon$) when the iteration terminated, and $\sigma_{\min}(\cdot)$ is the minimum singular value. From Table \ref{alg-perfor}, although solutions are obtained for all these four cases, Algorithm 3.5 performs much better for the latter two cases in terms of iteration numbers and singularity of the perturbed controllability matrices. This is consistent with the PTSC property of the corresponding structured systems. In fact, from the PTSC analysis in Example \ref{exp1-con}, the SCRPs for the first two cases are by no means feasible, while generically feasible for the third case. On the other hand, only after the algorithm parameters are selected properly, can we tell the feasibility of the associated SCRPs solely from solutions returned by the SCRP algorithm.
			\begin{center}
				\begin{table*}
					\caption{{Performance of Algorithm 3.5 in \cite{khare2012computing} for the SCRPs on system (\ref{exampsys})}} \label{alg-perfor}
					\begin{tabular}{ c|c|c|c|c }
						\hline
						Perturbation structure & Number of iterations & Terminated tolerance & $\sigma_{\min}({\cal C}(A,b))$ & $\sigma_{\min}({\cal C}(A+\Delta A,b+\Delta b))$ \\
						$[\Delta \bar A_1,\Delta \bar b_1]$ (complex-field) &   150000 &   $1.1818\times 10^{-2}$ & 5.3401 & 0.3791 \\
						$[\Delta \bar A_1,\Delta \bar b_1]$ (real-field) &  15723 &  $5.3185\times 10^{-11}$ & 5.3401 & 5.3036 \\
						$[\Delta \bar A_2,\Delta  \bar b_2]$ (complex-field) &  14 &  $3.4031\times 10^{-12}$ & 5.3401 &  0.0036 \\
						$[\Delta A_2,\Delta b_2]$ (real-field) & 13 &  $9.2831\times 10^{-13}$ & 5.3401 & 0.0103 \\
						\hline
					\end{tabular}
				\end{table*}
			\end{center}
			\vspace{-0.4cm}
			
			Several remarks about the application of PTSC are worthwhile:
			
		i) Problem (\ref{SCR}) is non-convex; even checking its feasibility for a specific perturbation structure is not an easy task.\footnote{A recent work \cite{zhang2022real} shows that checking the feasibility of a variant of the SCRP is NP-hard, if the perturbation matrix $[\Delta A, \Delta B]$ has an affine parameterization structure.}  Towards this problem, the majority of existing algorithms are numerical and iterative, and assume the feasibility as a default premise \cite{khare2012computing,bianchin2016observability,johnson2018structured,zhang2022real}. 
			Unlike those numerical algorithms, the proposed PTSC framework (valid in the generic sense for checking the problem feasibility) has the following advantages: it is guaranteed to be verified in $O(|{\cal E}_{\bar F}|n^4)$ time for single-input systems by resorting to some graph-theoretic algorithms, and it is free from the numerical difficulty of rounding errors. 
			

ii) Taking into account more prior knowledge on the values of the original system will certainly lead to a more precise answer to the feasibility of the SCRP. However, fixing the unknown generic entries to some constants will inevitably lead to results with a mix of numerical properties of the nomimal systems and combinatorial properties of the perturbation structure. By contrast, in the PTSC framework, all results rest on the combinatorial properties of the original systems and perturbation structure.

iii) We stress that the PTSC framework may be preferable in the context of network systems. A potential application scenario is that, one wants to 
estimate the controllability robustness of a network topology against the perturbations of a subset of its edges. Usually the exact edge weights are hard to obtain, but the interconnection topology is easily accessible. In this scenario, the PTSC framework could provide a generic answer for this controllability robustness estimation problem.

iv) Finally, we note that the PTSC could be extended to take more algebraic dependence among the unknown generic entries and the perturbed ones into account (see Sec \ref{sec_gener}). However, it is expected that, finding computationally efficient criteria for PTSC would then be more involved as one would need to redefine the associated graphs. Our results for PTSC could be a starting point towards some more realistic but complicated scenarios.

			\section{Conclusion}
			This paper proposes a novel notion, namely PTSC, to characterize the generic property of controllability preservation for realizations of a structured system under structured perturbations. The perturbation structure can be arbitrary relative to the structure of the original systems. A decomposition-based necessary and sufficient condition is given for single-input systems to be PTSC, which can be verified efficiently. Some intuitive graph-theoretic conditions for PSSC are subsequently given. It is shown PTSC could provide some feasibility condition for SCRPs from a generic view. Our future research directions include, exploring necessary and sufficient conditions for PTSC of multi-input systems, studying PTSC in the real field, and extending PTSC to consider parameter dependence.
			
			\vspace{-0.2cm}
			\begin{appendix}
				\subsection{Proof of Theorem \ref{generic_PSC_multi}}
				Our proof relies on an elementary tool from the theory on system of polynomial equations. Let
				\begin{equation}\label{sys_poly}f_1(x_1,...,x_n)=0,\cdots,f_t(x_1,...,x_n)=0 \end{equation}
				be a system of polynomial equations, where each $f_i$ ($i=1,...,t$) is a polynomial in the variables $(x_1,...,x_n)$ with real or complex coefficients. For notation simplicity, (\ref{sys_poly}) is denoted by $(f_1,...,f_t)$. Equation (\ref{sys_poly}) is said to be unsolvable, if it has no solution for $(x_1,...,x_n)\in {\mathbb C}^n$. The following result, known as Hilbert's Nullstellensatz, gives an elementary certificate for solvability of $(f_1,...,f_t)$.
				
				\begin{lemma}[Hilbert's Nullstellensatz]\citep[Chap 7.4]{Sturmfels2002} \label{Hilbert-Null} The system of polynomial equation (\ref{sys_poly}) is unsolvable, if and only if there exist some polynomials $g_i$ ($i=1,...,t$) in the variables $(x_1,...,x_t)$ with complex coefficients, such that
					$$\sum \nolimits_{i=1}^t g_if_i=1.$$
				\end{lemma}
				
				The following result provides an upper bound on the total degree of the $g_i$ (i.e., the highest degree of a monomial in $g_i$, denoted by $\deg(g_i)$) in Lemma \ref{Hilbert-Null}. Any such an upper bound is called an effective Nullstellensatz in the literature.
				
				\begin{lemma}[Effective Nullstellenstatz]\cite{kollar1988sharp}\label{effective-null} Let $f_1,...,f_t$ be polynomials in $n\ge 2$ variables $(x_1,...,x_n)$, of total degree $d_1,...,d_t$. If there are polynomials $g_i$ in $(x_1,...,x_n)$ such that $\sum \nolimits_{i=1}^t g_if_i=1$, then they can be chosen such that $\deg(f_ig_i)\le \max(3,d)^{\min(n,t)}$, where $d$ is the maximum of the degrees of the $f_i$ ($i=1,...,t$).
				\end{lemma}
				
				{\bf Proof of Theorem \ref{generic_PSC_multi}}: Let $p_1,...,p_r$ be values that the nonzero entries of $[A,B]$ take, and $\bar p_1,...,\bar p_l$ be values that the perturbed entries of $[\Delta A, \Delta B]$ take with $[\Delta A, \Delta B]\in {\bf S}_{\bar F}$, where $r$, $l$ are the number of nonzero entries in $[\bar A,\bar B]$ and $\bar F$, respectively. Let $p\doteq (p_1,...,p_r)$ and $\bar p\doteq (\bar p_1,...,\bar p_l)$. $(A+\Delta A, B+ \Delta B)$ is uncontrollable, if and only if the $n\times (nm)$ matrix $M(p,\bar p)\doteq {\cal C}(A+\Delta A, B+ \Delta B)$ is row rank deficient.  This is equivalent to that
				\begin{equation}\label{sys-poly-eq} \det(M(p,\bar p)[{\cal J}_n, {\cal I}])=0, \forall {\cal I}\in {\cal J}_{nm}^n,\end{equation}
				which induces a system of $N\doteq |{\cal J}_{nm}^n|={\tiny \left(\begin{array}{l} mn \\ n \end{array}\right)}$ polynomial equations in $\bar p$. For the $i$th element $\cal I$ in ${\cal J}_{nm}^n$, define $f_i\doteq \det(M(p,\bar p)[{\cal J}_n,{\cal I}])$. According to Lemma \ref{Hilbert-Null}, the system of polynomial equations (\ref{sys-poly-eq}) is unsolvable, if and only if there exist polynomials $g_i$ in $\bar p$ ($i=1,...,N$), such that
				\begin{equation} \label{null-formula}\sum \nolimits_{i=1}^N g_if_i=1.\end{equation}
				By Lemma \ref{effective-null},  such $g_i$ can be chosen so that
				$\deg(g_if_i)\le \max(3,n^2)^{\min(N,l)},$ where $n^2$ is the upper bound of $\deg(f_i)$. Therefore, $\deg(g_i)\le \max(3,n^2)^{\min(N,l)}-\deg(f_i)\doteq q_i$. Built on this, suppose there are coefficients $\{\alpha_{r_1^i\cdots r_l^i}\}$ with $(r_1^i,...,r_l^i)\in \Theta_i\doteq \left\{(r_1^i,...,r_l^i)\in {\mathbb N}^l: r_1^i+\cdots+r_l^i\le q_i\right\}$, such that $g_i$ can be expressed as
				\begin{equation}\label{expression-expand} g_i=\sum \nolimits_{(r_1^i,\cdots,r_l^i)\in \Theta_i} \alpha_{r_1^i\cdots r_l^i}\bar p_1^{r_1^i}\cdots \bar p_{l}^{r_l^i}.  \end{equation}
				Substituting (\ref{expression-expand}) into (\ref{null-formula}) and letting the coefficient of $\bar p_1^{r_1^i}\cdots \bar p_{l}^{r_l^i}$ for each distinct nonzero $(r_1^i,\cdots,r_l^i)\in \Theta_i$ equal zero, we obtain a finite system of linear equations in the variables $\{\alpha_{r_1^i\cdots r_l^i}\}$, which can be collectively expressed as
				\begin{equation} \label{eq-syslinear}
					H(p)\alpha=\left[
					\begin{array}{c}
						0_{(L-1)\times 1} \\
						1 \\
					\end{array}
					\right].
				\end{equation}
				In the formula above, $\alpha$ is a vector by stacking all $\{\alpha_{r_1^i\cdots r_l^i}|^{i=1,..., N}_{(r_1^i,\cdots,r_l^i)\in \Theta_i}\}$, $L$ is the total number of distinct elements in $\left\{\bar p_1^{r_1^i}\cdots \bar p_l^{r_l^i}: (r_1^i,\cdots,r_l^i)\in \Theta_i, i=1,...,N \right\}$, and $H(p)$ is a matrix whose entries are polynomials in $p$. Denote the right-hand side of (\ref{eq-syslinear}) as $\beta$. From the linear algebra theory \cite{Sturmfels2002}, Equation (\ref{eq-syslinear}) has a solution for $\alpha$, if and only if
				\begin{equation}\label{final-rank} {\rm rank}(H(p))-{\rm rank}([H(p),\beta])=0.\end{equation} 
				Note that both $H(p)$ and $[H(p),\beta]$ are polynomial-valued matrices of variables $p$. From the definition of generic rank, it holds that
				$$\begin{array}{l} {\rm rank}(H(p))= {\rm grank}(H(p)), \forall p\in {\mathbb C}^r\backslash {\cal P}_1, \\
					{\rm rank}([H(p),\beta])= {\rm grank}([H(p),\beta]), \forall p\in {\mathbb C}^r\backslash {\cal P}_2, \end{array}$$
				where ${\cal P}_1$ and ${\cal P}_2$ are proper varieties in ${\mathbb C}^r$. Hence, we have
				$$\begin{array}{l}{\rm rank}(H(p))-{\rm rank}([H(p),\beta])\\ = {\rm grank}(H(p))-{\rm grank}([H(p),\beta]), \forall p\in {\mathbb C}^r\backslash \left({\cal P}_1\cup {\cal P}_2\right).\end{array}$$
				Notice that ${\rm grank}(H(p))$ and ${\rm grank}([H(p),\beta])$ are independent of the values of $p$, and ${\cal P}_1\cup {\cal P}_2$ is still a proper variety of ${\mathbb C}^r$. Therefore, the value of the left-hand side of (\ref{final-rank}) is generic, which achieves the same number (either $0$ or $1$) for all $p\in {\mathbb C}^r$ except $p\in {\cal P}_1\cup {\cal P}_2$, and so is the solvability of Equation (\ref{eq-syslinear}). This means, the existence of $g_i$ for Equation (\ref{null-formula}) is generic, too. By Lemma \ref{Hilbert-Null}, solvability of the system of polynomial equations (\ref{sys-poly-eq}) is also generic for $p\in {\mathbb C}^r$. As $(\bar A,\bar B)$ is structurally controllable, the complement of the set $\{p\in {\mathbb C}^r: (A,B)\in {\bf CS}(\bar A, \bar B)\}$ has zero measure in ${\mathbb C}^r$. This immediately leads to the required statement in Theorem \ref{generic_PSC_multi}.  $\hfill\square$
				
				\subsection{Proof of Proposition \ref{mainpro}}  \label{mainpro-proof}
				\begin{lemma}\citep[Lem 2]{Rational_function} \label{root-factor} Let $p_1(\lambda,t_1,...,t_r)$ and $p_2(\lambda,t_1,...,t_r)$ be two polynomials on the variables $\lambda,t_1,...,t_r$ with real coefficients. Then, 1) For all $(t_1,...,t_r)\in {\mathbb C}^r$, $p_1(\lambda,t_1,...,t_r)$ and $p_2(\lambda,t_1,...,t_r)$ share a common zero for $\lambda$, if and only if $p_1(\lambda,t_1,...,t_r)$ and $p_2(\lambda,t_1,...,t_r)$ share a common factor in which the leading degree for $\lambda$ is nonzero; 2) If the condition above is not satisfied, then for almost al $(t_1,...,t_r)\in {\mathbb C}^r$ (except for a set with zero measure), $p_1(\lambda,t_1,...,t_r)$ and $p_2(\lambda,t_1,...,t_r)$ do not share a common zero for $\lambda$.
				\end{lemma}

				\begin{lemma} \label{matrix_pencil}
					Let $M$ be an $n\times n$ generic matrix over the variables $t_1,...,t_r$, and $E\in \{0,1\}^{n\times n}$, where each row, as well as each column of $E$, has at most one entry being $1$ and the rest being $0$. Let $P_\lambda\doteq M-\lambda E$ be a generic matrix pencil. Moreover, ${\cal B}({P_{\lambda}})$ is the bipartite graph associated with $P_{\lambda}$ defined similarly to ${\cal B}(H_\lambda)$ by replacing $I$ with $E$ (notably, self-loops are edges in ${\cal E}_E\cap {\cal E}_M$). The following are true:

					1) Suppose ${\cal B}({P_{\lambda}})$ contains no self-loop. Let $\gamma_{\min}({\cal B}({P_{\lambda}}))$ and $\gamma_{\max}({\cal B}({P_{\lambda}}))$ be respectively the minimum number of $\lambda$-edges and maximum number of $\lambda$-edges contained in a matching over all maximum matchings of ${\cal B}({P_{\lambda}})$. Then, the generic number (i.e., true for almost all values of $t_1,...,t_r$) of nonzero roots of $\det (P_{\lambda})$ for $\lambda$ (counting multiplicities) equals $\gamma_{\max}({\cal B}({P_{\lambda}}))- \gamma_{\min}({\cal B}({P_{\lambda}}))$.

					2) If ${\cal B}({P_{\lambda}})$ is DM-irreducible, then $\det (P_{\lambda})$ generically has nonzero roots for $\lambda$ whenever ${\cal B}({P_{\lambda}})$ contains a self-loop.
					
					3) Suppose that ${\cal B}({P_{\lambda}})$ is DM-irreducible. Let ${\cal T}_i$ be the subset of variables of $t_1,...,t_r$ that appear in the $i$th column of $M$. Then, for each $i\in \{1,...,n\}$, every nonzero root of $\det (P_{\lambda})$ (if exists) cannot be independent of ${\cal T}_i$.
				\end{lemma}
				
				\begin{proof}  We first prove an useful observation, that every maximum matching of ${\cal B}({P_{\lambda}})$ corresponds to a nonzero term (monomial) in $\det (P_{\lambda})$ that cannot be zeroed out by other terms, which is crucial to the following proofs. For this purpose, consider a term $\lambda^{r_1}t_{l_1}t_{l_2}\cdots t_{l_{r_2}}$ associated with a maximum matching of ${\cal B}({P_{\lambda}})$, where $r_1+r_2=n$ and $\{l_1,...,l_{r_2}\}\subseteq \{1,...,r\}$. The only case to zero out $\lambda^{r_1} t_{l_1}t_{l_2}\cdots t_{l_{r_2}}$ in $\det (P_{\lambda})$ is that there exists a term being $\lambda^{r_1} t_{l_1}t_{l_2}\cdots t_{l_{r_2}}$ associated with another maximum matching of ${\cal B}({P_{\lambda}})$.  Now suppose that $R$ (resp. $C$) is the set of row indices (resp. column indices) of $t_{l_1},...,t_{l_{r_2}}$ in $M$, recalling that each $t_i$ ($i=1,...,r$) appears only once. Then, $\lambda^{r_1}$ must correspond to $r_1$ `1' entries in $E[{\cal J}_n\backslash R, {\cal J}_n\backslash C]$. However, as each row and each column has at most one `$1$' in $E$, the aforementioned configuration for `$1$' entries is unique, contracting the existence of two different maximum matchings associated with $\lambda^{r_1} t_{l_1}t_{l_2}\cdots t_{l_{r_2}}$.
					
					We now prove 1). Suppose $\det (P_{\lambda})$ can be factored as $\lambda^{l}f(\lambda,t_1,...,t_r)$, where $l\in {\mathbb N}$ and $f(\lambda,t_1,...,t_r)$ is a polynomial of $\lambda,t_1,...,t_r$ that does not contain factors in the form of $\lambda^{\bar l}$ for any $\bar l\ge 1$. Because of the above observation, every term associated with a maximum matching of ${\cal B}({P_{\lambda}})$ must contain the factor $\lambda^l$. Therefore, the number $l$ of zero roots of $\det (P_{\lambda})$ equals  $\gamma_{\min}({\cal B}({P_{\lambda}}))$.  In addition, the maximum degree of $\lambda$ in $\lambda^{l}f(\lambda,t_1,...,t_r)$ appears in a term associated with a maximum matching containing the maximum number of $\lambda$-edges, which is exactly $\gamma_{\max}({\cal B}({P_{\lambda}}))$. The conclusion in 1) then follows from the fundamental theorem of algebra.
					
					Next, we prove 2). Consider a self-loop with the entry being $t_l-\lambda$ ($1\le l \le r$). As ${\cal B}({P_{\lambda}})$ is DM-irreducible, every nonzero entry must be contained in $\det (P_{\lambda})$ by Definition \ref{DM-def}, which means $t_l-\lambda$ is contained in some term $(t_l-\lambda)f$ of $\det (P_{\lambda})$, where $f$ denotes a polynomial over variables $\{t_1,...,t_r\}\backslash \{t_l\}$ and $\lambda$. This term can be written as the sum of two terms $t_lf$ and $-\lambda f$, which indicates $\det (P_{\lambda})$ contains at least two monomials whose degrees for $\lambda$ differ from each other. Then, following the similar reasoning to the proof of 1),  $\det (P_{\lambda})$ contains at least one nonzero root.
					
					We are now proving 3). Suppose such a nonzero root exists that is independent of ${\cal T}_i$ for some $i\in \{1,...,n\}$, and denote it by $z$. Let ${\cal T}[{\cal I}_1,{\cal I}_2]$ be the set of variables in $t_1,...,t_r$ that appear in $M[{\cal I}_1,{\cal I}_2]$ for ${\cal I}_1,{\cal I}_2\subseteq {\cal J}_n$, and let $R({\cal T}_s)$ (resp. $C({\cal T}_s)$) be the set of row (resp. column) indices of variables ${\cal T}_s\subseteq \{t_1,...,t_r\}$.  Suppose $[P_{\lambda}]_{k_0,i}=\lambda$ for some $k_0\in {\cal J}_n\backslash R({\cal T}_i)$ ($k_0$ can be empty). Upon letting all $t_k\in {\cal T}_i$ be zero, we obtain ($P_z=M-z E$)
					$${\small\begin{array}{c}\begin{aligned} &\det (P_z)=\sum\nolimits_{j=1}^n(-1)^{i+j}[P_{z}]_{ji}\det(P_z[{\cal J}_n\backslash \{j\}, {\cal J}_n\backslash \{i\}])\\
								&=z\cdot \det (P_z[{\cal J}_n\backslash \{k_0\}, {\cal J}_n\backslash \{i\}])=0, \end{aligned}\end{array}}$$
					which indicates
					\begin{equation}\label{lambda_zero} \det (P_z[{\cal J}_n\backslash \{k_0\}, {\cal J}_n\backslash \{i\}])=0,\end{equation}
					as $z\ne 0$. Since $P_{\lambda}[{\cal J}_n\backslash \{k_0\}, {\cal J}_n\backslash \{i\}]$ has full generic rank from Lemma \ref{reduciblility}\footnote{Hereafter, when referring to the generic rank of $P_\lambda$ (or its submatrices), we should regard $P_\lambda$ as a polynomial-valued matrix of the free values $t_1,...,t_r$ and $\lambda$.}, it concludes that $z$ depends solely on the variables ${\cal T}[{\cal J}_n\backslash \{k_0\}, {\cal J}_n\backslash \{i\}]$. Similarly, because of (\ref{lambda_zero}), for each $t_k\in {\cal T}_i$, fixing all $t_j\in {\cal T}_i\backslash \{t_k\}$ to be zero yields
					$$\det (P_z[{\cal J}_n\backslash R(\{t_k\}), {\cal J}_n\backslash \{i\}])=0,$$
					which indicates that $z$ depends on the variables ${\cal T}[{\cal J}_n\backslash R(\{t_k\}), {\cal J}_n\backslash \{i\}]$, being independent of the remaining variables. Taking the intersection of ${\cal T}[{\cal J}_n\backslash \{j\}, {\cal J}_n\backslash \{i\}]$ over all $j\in R({\cal T}_i)\cup \{k_0\}$, we obtain ${\cal T}[\Theta,{\cal J}_n\backslash \{i\}]$, where $\Theta\doteq {\cal J}_n\backslash (R({\cal T}_i)\cup \{k_0\})$ (see Fig. \ref{proof-illustration} for illustration). That is, $z$ depends on variables ${\cal T}[\Theta,{\cal J}_n\backslash \{i\}]$, and makes $P_{\lambda}[\Theta,{\cal J}_n\backslash \{i\}]$ row rank deficient. However, for each pair $(j,l)$, $j\! \in\! {\cal J}_n\backslash \{i\}$, $l\in  R({\cal T}_i)$, it~holds
					$$\begin{array}{c}\begin{aligned} &{\rm grank} (P_{\lambda}[\Theta,{\cal J}_n\backslash \{i,j\}])\\
							&\mge {\rm grank}((M\!-\!\lambda E)[{\cal J}_n\backslash \{l\} ,{\cal J}_n\backslash \{j\}])\!-\!(|R({\cal T}_i)\cup\{k_0\}|\!-1\!)\\
							&\myeq n-|R({\cal T}_i)\cup\{k_0\}|=|\Theta|,\end{aligned}\end{array}$$
					where (a) is due to that $P_{\lambda}[\Theta,{\cal J}_n\backslash \{i,j\}]$ is obtained by deleting $|R({\cal T}_i)\cup\{k_0\}|\!-1\!$ rows from $P_\lambda[{\cal J}_n\backslash \{l\} ,{\cal J}_n\backslash \{j\}]$ (cf. Fig. \ref{proof-illustration}), and (b) comes from ${\rm grank}(P_{\lambda}[{\cal J}_n\backslash \{l\} ,{\cal J}_n\backslash \{j\}])=n-1$ by 2) of Lemma \ref{reduciblility}. That is, after deleting any column from $P_{\lambda}[\Theta,{\cal J}_n\backslash \{i\}]$, the resulting matrix remains of full row generic rank, which induces at least one nonzero polynomial equation constraint on $z$ and ${\cal T}[\Theta,{\cal J}_n\backslash \{i,j\}]$. This indicates $z$ depends on ${\cal T}[\Theta,{\cal J}_n\backslash \{i,j\}]$, or equivalently, being independent of ${\cal T}[\Theta,\{j\}]$, for each $j\in {\cal J}_n\backslash \{i\}$. It finally concludes that $z$ is independent of the variables ${\cal T}[\Theta, {\cal J}_n\backslash \{i\}]$, causing a contraction. Therefore, the assumed $z$ cannot exist.
				\end{proof}
				
				\begin{figure}
					\centering
					\includegraphics[width=1.5in]{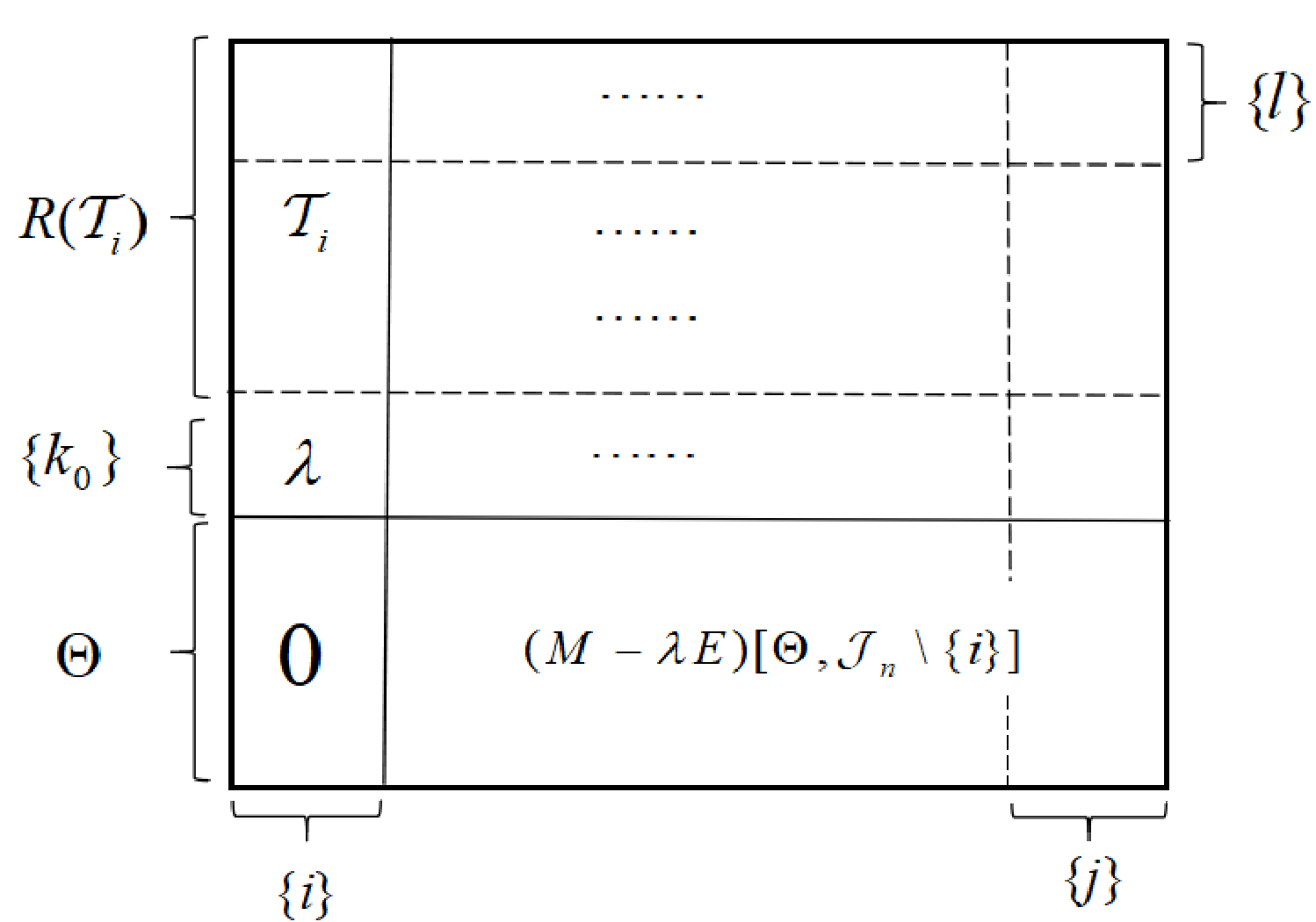}\\
					\caption{Illustration of proof for property 3) of Lemma \ref{matrix_pencil}}\label{proof-illustration}
				\end{figure}
				
				\begin{lemma}\label{immediate} Let $M^{j_{\rm c}}_\lambda$ and ${\Omega_j}$ be defined in (\ref{DMdecomposition}) and (\ref{important-set}). For each $k\in {\Omega_j}$, let $\tilde{M}^{j_{\rm c}}_{ki^*}(\lambda)\doteq M_{\lambda}^{j_{\rm c}}[\bar {\cal V}^+_{ki^*},\bar {\cal V}^-_{ki^*}]$, i.e.,
					{\small$$
						\tilde{M}^{j_{\rm c}}_{ki^*}(\lambda)\!=\!\left[\!
						\begin{array}{ccc}
							M_k^{j_{\rm c}}(\lambda) & \cdots & M_{ki^*}^{j_{\rm c}}(\lambda) \\
							0 & \ddots & \vdots \\
							0 & 0 & M^{j_{\rm c}}_{i^*}(\lambda)[{\cal J}_{|{\cal V}^+_{i^*}|}\backslash \{\bar i\},{\cal J}_{|{\cal V}^+_{i^*}|}]
						\end{array}
						\!\right].$$}Then, $\tilde{M}^{j_{\rm c}}_{ki^*}(\lambda)$ generically has full row rank when $\lambda \in \{z\in {\mathbb C}\backslash \{0\}: \det (M_k^{j_{\rm c}}(z))=0 \}$, if and only if ${\rm minWMM}({\cal G}^{j_{\rm c}}_{ki^*})< |{\cal V}^{+}_k|$.
				\end{lemma}

				
				\begin{proof}  Sufficiency: Let $n_1=|{\cal V}^+_k \cup {\cal V}^+_{k+1}\cdots \cup{\cal V}^+_{i^*}|$. By Lemma \ref{reduciblility}, ${\rm mt}({\cal G}^{j_c}_{i^*}-\{x_i\})=|{\cal V}^+_{i^*}|-1$. Consequently, ${\cal G}^{j_{\rm c}}_{ki^*}$ has a maximum matching with size $n_1-1$ from its structure. Suppose that ${\cal G}^{j_{\rm c}}_{ki^*}$ has a maximum matching with weight less than $|{\cal V}^{+}_k|$. Then, ${\cal B}(\tilde{M}^{j_{\rm c}}_{ki^*}(\lambda)[{\cal J}_{n_1-1},{\cal J}_{n_1}\backslash \{\bar k\}])$ for some $\bar k\in \{1,...,|{\cal V}^-_k|\}$ must have a matching with size $n_1-1$. Indeed, if this is not true, then any maximum matching of ${\cal G}^{j_{\rm c}}_{ki^*}$ must matches ${\cal V}^-_k$, which certainly leads to a weight equaling $|{\cal V}^+_k|$, noting that each edge not incident to ${\cal V}^-_k$ has a zero weight. Furthermore, due to the DM-irreducibility of ${\cal G}^{j_{\rm c}}_k$, from Lemma \ref{matrix_pencil}, any nonzero root of $\det (M_k^{j_{\rm c}}(\lambda))$ cannot be independent of the variables in the $\bar k$th column of $M_k^{j_{\rm c}}(\lambda)$.\footnote{Note the case where $|{\cal V}^+_k|=1$ and $M_k^{j_{\rm c}}(\lambda)=\lambda$ has been excluded by the nonzero root assumption.} Therefore, $\det (M_k^{j_{\rm c}}(\lambda))$ and $\det (\tilde{M}^{j_{\rm c}}_{ki^*}(\lambda)[{\cal J}_{n_1-1},{\cal J}_{n_1}\backslash \{\bar k\}])$ generically do not share a common nonzero root, since the latter determinant cannot contain the variables in the $\bar k$th column of $M_k^{j_{\rm c}}(\lambda)$ (except for $\lambda$).  That is, $\det (\tilde{M}^{j_{\rm c}}_{ki^*}(\lambda)[{\cal J}_{n_1-1},{\cal J}_{n_1}\backslash \{\bar k\}])$ is generically nonzero for $\lambda \in \{z\in {\mathbb C}\backslash \{0\}: \det (M_k^{j_{\rm c}}(z))=0\}$, leading to the full row rank of $\tilde{M}^{j_{\rm c}}_{ki^*}(\lambda)$.
					
					Necessity: If $k=i^*$, the necessity is obvious. Consider $k<i^*$. Suppose ${\rm minWMM}({\cal G}^{j_{\rm c}}_{ki^*})= |{\cal V}^{+}_k|$. Then,  based on the above analysis, any maximum matching of ${\cal G}^{j_{\rm c}}_{ki^*}$ must match ${\cal V}^-_k$, which leads to a zero determinant of any $(n_1-1)\times (n_1-1)$ submatrix of ${\tilde M}^{j_{\rm c}}_{ki^*}(\lambda)$, due to the block-triangular structure of ${\tilde M}^{j_{\rm c}}_{ki^*}(\lambda)$ and the fact that $\det (M_k^{j_{\rm c}}(\lambda))=0$, contradicting the full row rank of ${\tilde M}^{j_{\rm c}}_{ki^*}(\lambda)$.
				\end{proof}

				{\bf Proof of Proposition \ref{mainpro}}: In the following, suppose
				vertex $x_i$ corresponds to the $\hat i$th row of $M^{j_{\rm c}}_\lambda$ after the row permutation by $P$, i.e. $[M_\lambda^j]_{\hat i}=[PH_\lambda^j]_{\hat i}=[H_\lambda^j]_i$. Recall the involved $H_\lambda, M^{j_c}_{\lambda}$ and their submatrices are treated as generic matrix pencils.
				
				Sufficiency: From Lemma \ref{matrix_pencil}, we know that for each $k\in \Omega_j$, there exists a nonzero $\lambda$ making $\det (M_k^{j_{\rm c}}(\lambda))=0$.  To distinguish such value from the variable $\lambda$, we denote it by $z$ (i.e., $z\ne 0$ and $\det (M_k^{j_{\rm c}}(z))=0$). From Lemma \ref{root-factor}, it generically holds that $\det (M_l^j(z))\ne 0$ for all $l\in \{1,...,d\}\backslash \{k\}$, as $\det (M_k^{j_{\rm c}}(\lambda))$ and $\det (M_l^j(\lambda))$ do not share any common factor except the power of $\lambda$. Due to the block-triangular structure of $M^{j_{\rm c}}_z$ (obtained by replacing $\lambda$ with $z$ in $M^{j_{\rm c}}_\lambda$), it can be seen readily that if $\tilde{M}^{j_{\rm c}}_{ki^*}(z)$ defined in Lemma \ref{immediate} generically has full row rank, then $M^{j_{\rm c}}_z[{\cal J}_n\backslash \{\hat i\}, {\cal J}_{n}]$ will do. The former condition has been proven in Lemma \ref{immediate}.
				
				Note also that $M^{j_{\rm c}}_z$ generically has rank $n-1$ as otherwise $[A-zI,b]$ generically has rank less than $n$, contradicting the structural controllability of $(\bar A, \bar b)$. Therefore, from Lemma \ref{basiclemma}, letting ${\hat q}$ be a nonzero vector in the left null space of $M^{j_{\rm c}}_z$, we have ${\hat q}_{\hat i}\ne 0$. For the generic realization $[A,b]$ of $[\bar A, \bar b]$, by letting $[\Delta A, \Delta b]_{ij}=-1/{\hat q}_{\hat i} \sum \nolimits_{l=1}^n{\hat q}_l(P[A-zI,b])_{lj},$ we get
				$$\begin{array}{c}\begin{aligned}&{\hat q}^{\intercal}P([A-zI,b]+[\Delta A, \Delta b])[{\cal J}_n,\{j\}]\\
						&={\hat q}_{\hat i}(P[\Delta A, \Delta b])_{\hat i,j}+\sum \nolimits_{l=1}^n {\hat q}_l(P[A-zI,b])_{lj}\\
						&=0,\end{aligned} \end{array}$$
				where the second equality is due to $(P[\Delta A,\Delta b])_{\hat i,j}=[\Delta A, \Delta b]_{ij}$. Upon defining $q^{\intercal}\doteq {\hat q}^{\intercal}P$, we have
				$q^{\intercal}([A-zI,b]+[\Delta A, \Delta b])=0,$ which comes from the fact  $q^{\intercal}H_z^{j_{\rm c}}Q=0$ and $Q$ is invertible.
				
				Necessity: For the existence of $q$ making the condition in Proposition \ref{mainpro} satisfied, it is necessary $H_\lambda^{j_{\rm c}}$ should be of rank deficient at some nonzero value for $\lambda$ generically. Denote such a value by $z$ for the sake of distinguishing it from the variable $\lambda$. Since DM-decomposition does not alter the rank, $M_{z}^{j_{\rm c}}$ should be of row rank deficient too. From the block-triangular structure of $M_z^{j_{\rm c}}$, there must exist some $k\in \{1,...,d\}$, such that $M^{j_{\rm c}}_{k}(z)$ is singular generically. From Lemma \ref{matrix_pencil}, such an integer $k$ must correspond to a ${\cal G}_k^{j_{\rm c}}$ satisfying ${\gamma}_{\rm nz}({\cal G}_k^{j_{\rm c}})=1$. We consider two cases: i) $k> i^*$, and ii) $k\le i^*$.
				
				In case i), since $k>i^*$, from the upper block-triangular structure of $M_z^{j_{\rm c}}$, it is clear that $M_z^{j_{\rm c}}[{\cal J}_n\backslash \{\hat i\}, {\cal J}_n]$ is of row rank deficient when $\det (M^{j_{\rm c}}_k(z))=0$. Note that it holds generically ${\rm rank}(M_{z}^{j_{\rm c}})\ge n-1$, as otherwise ${\rm rank}(H_z)< n$, which is contradictory to the structural controllability of $(\bar A, \bar b)$. Consequently, $M_{z}^{j_{\rm c}}$ has a left null space with dimension one. Denote by $\hat q$ the vector spanning that space. From Lemma \ref{basiclemma}, ${\hat q}_{\hat i}=0$. As a result, for any $[\Delta A, \Delta  b]\in {\bf S}_{[ {\Delta \bar A},  {\Delta \bar b}]}$,
				$$\hat q^{\intercal} \big\{M^j_{z}+(P[\Delta A, \Delta  b])[{\cal J}_n,\{j\}]\big\}\myeqa \hat q^{\intercal} M^j_z \ne 0,$$
				where (a) results from $\hat q^\intercal_{\hat i}(P[\Delta A, \Delta b])_{\hat ij}=0$, and the inequality  from $\hat q^{\intercal} [M^{j_{\rm c}}_z,M^j_z]\ne 0$, as otherwise $\hat q^{\intercal} P H_z=0$ meaning that $z$ will be an uncontrollable mode. Hence, case i) cannot lead to the required results. 
				
				Therefore, $k$ must fall into case ii). Now suppose ${\rm minWMM}({\cal G}^{j_{\rm c}}_{ki^*})= |{\cal V}^{+}_k|$. Then, from Lemma \ref{immediate} and by the block-triangular structure of $M_z^{j_{\rm c}}$, we obtain that $M_z^{j_{\rm c}}[{\cal J}_n\backslash \{\hat i\}, {\cal J}_n]$ is generically of row rank deficient. By Lemma \ref{basiclemma} and following the similar reasoning to case i), it turns out that the requirement in Proposition \ref{mainpro} cannot be satisfied. This proves the necessity. $\hfill\square$
			\end{appendix}
	}}
	
	{
		\bibliographystyle{ieeetr}
		{\footnotesize
			\bibliography{yuanz3}

\begin{thebibliography}{10}

\bibitem{Buldyrev2009Catastrophic}
S.~V. Buldyrev, P.~Roni, P.~Gerald, S.~H~Eugene, and H.~Shlomo, ``Catastrophic
  cascade of failures in interdependent networks,'' {\em Nature}, vol.~464,
  no.~7291, pp.~1025--8, 2009.

\bibitem{fawzi2014secure}
H.~Fawzi, P.~Tabuada, and S.~Diggavi, ``Secure estimation and control for
  cyber-physical systems under adversarial attacks,'' {\em IEEE Transactions on
  Automatic Control}, vol.~59, no.~6, pp.~1454--1467, 2014.

\bibitem{mitra2019byzantine}
A.~Mitra and S.~Sundaram, ``Byzantine-resilient distributed observers for lti
  systems,'' {\em Automatica}, vol.~108, p.~108487, 2019.

\bibitem{zhang2020generic}
Y.~Zhang, Y.~Xia, J.~Zhang, and J.~Shang, ``Generic detectability and
  isolability of topology failures in networked linear systems,'' {\em IEEE
  Transactions on Control of Network Systems}, vol.~8, no.~1, pp.~500 -- 512,
  2021.

\bibitem{FabioFragility2018}
P.~Fabio, C.~Favaretto, S.~Zhao, and S.~Zampieri, ``Fragility and
  controllability tradeoff in complex networks,'' in {\em American Control
  Conference}, pp.~216--227, IEEE, 2018.

\bibitem{De2015Input}
C.~De~Persis and P.~Tesi, ``Input-to-state stabilizing control under
  denial-of-service,'' {\em IEEE Transactions on Automatic Control}, vol.~60,
  no.~11, pp.~2930--2944, 2015.

\bibitem{commault2008observability}
C.~Commault, J.-M. Dion, and D.~H. Trinh, ``Observability preservation under
  sensor failure,'' {\em IEEE Transactions on Automatic Control}, vol.~53,
  no.~6, pp.~1554--1559, 2008.

\bibitem{Rahimian2013Structural}
M.~A. Rahimian and A.~G. Aghdam, ``Structural controllability of multi-agent
  networks: Robustness against simultaneous failures,'' {\em Automatica},
  vol.~49, no.~11, pp.~3149--3157, 2013.

\bibitem{zhang2019minimal}
Y.~Zhang and T.~Zhou, ``Minimal structural perturbations for controllability of
  a networked system: Complexities and approximations,'' {\em International
  Journal of Robust and Nonlinear Control}, vol.~29, no.~12, pp.~4191--4208,
  2019.

\bibitem{generic}
J.~M. Dion, C.~Commault, and J.~Van~DerWoude, ``Generic properties and control
  of linear structured systems: {a} survey,'' {\em Automatica}, vol.~39,
  pp.~1125--1144, 2003.

\bibitem{paige1981properties}
C.~Paige, ``Properties of numerical algorithms related to computing
  controllability,'' {\em IEEE Transactions on Automatic Control}, vol.~26,
  no.~1, pp.~130--138, 1981.

\bibitem{WickDistance1991}
M.~Wicks and R.~DeCarlo, ``Computing the distance to an uncontrollable
  system,'' {\em IEEE Transactions on Automatic Control}, vol.~36, no.~1,
  pp.~39--49, 1991.

\bibitem{Hu2004Real}
G.~Hu and E.~J. Davison, ``Real controllability/stabilizability radius of lti
  systems,'' {\em IEEE Transactions on Automatic Control}, vol.~49, no.~2,
  pp.~254--257, 2004.

\bibitem{khare2012computing}
S.~R. Khare, H.~K. Pillai, and M.~N. Belur, ``Computing the radius of
  controllability for state space systems,'' {\em Systems \& Control Letters},
  vol.~61, no.~2, pp.~327--333, 2012.

\bibitem{johnson2018structured}
S.~C. Johnson, M.~Wicks, M.~{\v{Z}}efran, and R.~A. DeCarlo, ``The structured
  distance to the nearest system without property {P},'' {\em IEEE Transactions
  on Automatic Control}, vol.~63, no.~9, pp.~2960--2975, 2018.

\bibitem{bianchin2016observability}
G.~Bianchin, P.~Frasca, A.~Gasparri, and F.~Pasqualetti, ``The observability
  radius of networks,'' {\em IEEE transactions on Automatic Control}, vol.~62,
  no.~6, pp.~3006--3013, 2016.

\bibitem{zhang2022real}
Y.~Zhang, Y.~Xia, and Y.~Zhan, ``On real structured
  controllability/stabilizability/stability radius: Complexity and unified
  rank-relaxation based methods,'' {\em arXiv preprint arXiv:2201.01112}, 2022.

\bibitem{mayeda1979strong}
H.~Mayeda and T.~Yamada, ``Strong structural controllability,'' {\em SIAM
  Journal on Control and Optimization}, vol.~17, no.~1, pp.~123--138, 1979.

\bibitem{bowden2012strong}
C.~Bowden, W.~Holderbaum, and V.~M. Becerra, ``Strong structural
  controllability and the multilink inverted pendulum,'' {\em IEEE transactions
  on automatic control}, vol.~57, no.~11, pp.~2891--2896, 2012.

\bibitem{monshizadeh2014zero}
N.~Monshizadeh, S.~Zhang, and M.~K. Camlibel, ``Zero forcing sets and
  controllability of dynamical systems defined on graphs,'' {\em IEEE
  Transactions on Automatic Control}, vol.~59, no.~9, pp.~2562--2567, 2014.

\bibitem{mousavi2017structural}
S.~S. Mousavi, M.~Haeri, and M.~Mesbahi, ``On the structural and strong
  structural controllability of undirected networks,'' {\em IEEE Transactions
  on Automatic Control}, vol.~63, no.~7, pp.~2234--2241, 2017.

\bibitem{popli2019selective}
N.~Popli, S.~Pequito, S.~Kar, A.~P. Aguiar, and M.~Ili{\'c}, ``Selective strong
  structural minimum-cost resilient co-design for regular descriptor linear
  systems,'' {\em Automatica}, vol.~102, pp.~80--85, 2019.

\bibitem{jia2020unifying}
J.~Jia, H.~J. van Waarde, H.~L. Trentelman, and M.~K. Camlibel, ``A unifying
  framework for strong structural controllability,'' {\em IEEE Transactions on
  Automatic Control}, vol.~66, no.~1, pp.~391--398, 2020.

\bibitem{zhang2021ptsc}
Y.~Zhang and Y.~Xia, ``{PTSC}: a new notion for structural controllability
  under structured perturbations,'' in {\em 2021 40th Chinese Control
  Conference (CCC)}, pp.~4919--4924, IEEE, 2021.

\bibitem{Murota_Book}
K.~Murota, {\em Matrices and Matroids for Systems Analysis}.
\newblock Springer Science Business Media, 2009.

\bibitem{dummit2004abstract}
D.~S. Dummit and R.~M. Foote, {\em Abstract algebra}, vol.~3.
\newblock Wiley Hoboken, 2004.

\bibitem{C.T.1974Structural}
C.~T. Lin, ``Structural controllability,'' {\em IEEE Transactions on Automatic
  Control}, vol.~48, no.~3, pp.~201--208, 1974.

\bibitem{zhang2019structural}
Y.~Zhang and T.~Zhou, ``Structural controllability of an {NDS} with {LFT}
  parameterized subsystems,'' {\em IEEE Transactions on Automatic Control},
  vol.~64, no.~12, pp.~4920--4935, 2019.

\bibitem{kiraly2015algebraic}
F.~J. Kir{\'a}ly, L.~Theran, and R.~Tomioka, ``The algebraic combinatorial
  approach for low-rank matrix completion.,'' {\em J. Mach. Learn. Res.},
  vol.~16, no.~1, pp.~1391--1436, 2015.

\bibitem{zhang2021generic}
Y.~Zhang, Y.~Xia, H.~Zhang, G.~Wang, and L.~Dai, ``On the generic structured
  low-rank matrix completion,'' {\em arXiv preprint arXiv:2102.11490}, 2021.

\bibitem{Sturmfels2002}
B.~Sturmfels, {\em Solving Systems of Polynomial Equations}.
\newblock No.~97, American Mathematical Soc., 2002.

\bibitem{Geor1993Graph}
A.~George, J.~R. Gilbert, and J.~W.~H. Liu, {\em Graph Theory and Sparse Matrix
  Computation}.
\newblock Springer-Verlag: New York, 1993.

\bibitem{DB_West_graph}
D.~B. West, {\em Introduction to Graph Theory}.
\newblock Prentice hall, 2001.

\bibitem{liu2016control}
Y.-Y. Liu and A.-L. Barab{\'a}si, ``Control principles of complex systems,''
  {\em Reviews of Modern Physics}, vol.~88, no.~3, p.~035006, 2016.

\bibitem{KarowStructured2009}
M.~Karow and D.~Kressner, ``On the structured distance to uncontrollability,''
  {\em Systems \& Control Letters}, vol.~58, no.~2, pp.~128--132, 2009.

\bibitem{kollar1988sharp}
J.~Koll{\'a}r, ``Sharp effective nullstellensatz,'' {\em Journal of the
  American Mathematical Society}, pp.~963--975, 1988.

\bibitem{Rational_function}
K.~S. Lu and J.~N. Wei, ``Rational function matrices and structural
  controllability and observability,'' {\em IET Control Theory and
  Applications}, vol.~138, no.~4, pp.~388--394, 1991.

\end{thebibliography}
	}}

\end{document}